# PUZZLE GEOMETRY AND RIGIDITY: THE FIBONACCI CYCLE IS HYPERBOLIC

DANIEL SMANIA

ABSTRACT. We describe a new and robust method to prove rigidity results in complex dynamics. The new ingredient is the geometry of the critical puzzle pieces: under control of geometry and "complex bounds", two generalized polynomial-like maps which admits a topological conjugacy, quasiconformal outside the filled-in Julia set are, indeed, quasiconformally conjugated. The proof uses a new abstract removability-type result for quasiconformal maps, following ideas of Heinonen & Koskela and Kallunki & Koskela, optimized to complex dynamics. As the first application to this new method, we prove that, for even criticalities distinct of two, the period two cycle of the Fibonacci renormalization operator is hyperbolic with one-dimensional unstable manifold. To derive the exponential contraction in the hybrid classes, we use the non existence of invariant line fields in the Fibonacci tower, the topological convergence (both results by van Strien & Nowicki) and a new argument, distinct of the C. McMullen and M. Lyubich previous methods in the classic renormalization operator theory. We also describe other future applications.

## Contents

1991 *Mathematics Subject Classification.* 37F30, 37C15, 30C62, 30C65, 37F25 , 37F45, 37E20.
This work was supported by CNPq-Brazil grant 200764/01-2. Visits to Warwick University and ICTP-Trieste were supported by CNPq-Brazil grant 460110/00-4, Warwick University and ICTP-Trieste.

1. Introduction and results

1.1. **Complex dynamics results.** We will first introduce the complex dynamics results. If the reader is more interested in the removability-type result to quasiconformal maps, he/she must jump to section 1.2. Denote by $m(\cdot)$ the 2-dimensional Lebesgue measure.

**Definition 1.1.** *We say that $f \colon \cup_{i \leq n} U_i \to V$ is a **generalized polynomial like map** if*
- *$U_i$, $i \leq n$, and $V$ are simply connected domains on the complex plane. Moreover $U_i$ is compactly contained in $V$ and $m(\partial V) = m(\partial U_i) = 0$,*
- *$f$ has an analytic extension to a neighborhood of $\cup_i \overline{U}_i$.*
- *$f \colon U_0 \to V$ is a ramified covering with an unique critical point,*
- *For each $i \neq 0$, $f$ is an univalent map of $U_i$ to $V$,*

The simply connected components of the domain of $f^n$, for some $n > 0$, will be called **puzzle pieces at level** $n$ associated to $f$. If the critical point does not **escape** ($f^n(0)$ is well defined, for every $n \geq 0$), the **critical piece at level** $n$, denoted by $C_n$, is the connected component of the domain of $f^n$ which contains the critical point.

Rigidity problems appears in complex dynamics in many situations: the real Fatou conjecture to the quadratic family (Lyubich[Lyu97] and Graczyk & Swiatek[GS97]) and the hyperbolicity of the renormalization horseshoe (Lyubich[Lyu99]) were proved reducing these questions to a rigidity result. Usually the hard part is to conclude that two complex dynamical systems are quasiconformally conjugated, once they are topologically conjugated. The first step is to replace the topological conjugacy by another one which is quasiconformal outside the 'complicated' part of the dynamics (the Julia set), what, most of the time, is easy. Our first dynamical result is a method to carry out the next step:

**Theorem 1** (Rigidity theorem). *Let $f_j \colon \cup_i U_i^j \to V^j$, with $j = 1, 2$, be generalized polynomial like maps, whose critical point $0$ does not escape, such that*
- *There exists $\lambda < 1$ and an increasing sequence $n_k$ so that, if $C_n^j$ denotes the critical puzzle piece at level $n$ of $f_j$, then*
  - *diam $C_n^j \to_{n \to \infty} 0$,*
  - ***puzzle geometry control:*** *for any $j = 1, 2$, $k > 0$, there exists $x$ so that*
  $$B(x, \lambda \cdot \text{diam } C_{n_k}^j) \subset C_{n_k}^j$$
  - ***complex bounds:*** *for any $j = 1, 2$, $k > 0$, there exists an annulus $A_k^j \subset V^j$ such that*
    1. *mod $A_k^j \geq K_2$.*
    2. *The internal boundary of $A_k^j$ is $\partial C_{n_k}^j$,*

3. *For any $x$ which is not contained in $C_{n_k}^j$ and $m > 0$ minimal so that $f_j^m(x) \in C_{n_k}^j$, there exists a neighborhood $V_j(x,m)$ of $x$ so that*
   * $f_j^m(V_j(x,m)) = C_{n_k}^j \cup A_k^j$,
   * $f_j^m$ *is univalent in $V_j(x,m)$.*

Let $h \colon \mathbb{C} \to \mathbb{C}$ be a homeomorphism so that
- $h$ is quasiconformal in $\mathbb{C} - J(f_1)$,
- $h(\partial U_i^1) = \partial U_i^2$ and $h(\partial V^1) = \partial V^2$,
- $h$ is a conjugacy between $f_1$ and $f_2$.

Then $h \colon \mathbb{C} \to \mathbb{C}$ is quasiconformal.

**Remark 1.1.** *Indeed we can assume that the conjugacy $h$ is defined only in $\mathbb{C} - J(f_1)$, once the "complex bounds" condition on $f_1$ and $f_2$ implies that $h$ will have a continuous extension to $K(f_1) = J(f_1)$, which will be a conjugacy on $J(f_1)$.*

**Remark 1.2.** *The reader can state a similar theorem with many critical points involved.*

**Remark 1.3.** *We can obtain the third property of the annulus $A_k^j$ by the a stronger property, but useful in applications: $A_k^j \cap P(f) = \phi$ (called **unbranched complex bounds**). Here $P(f) := \overline{\{f^i(0) \colon i > 0\}}$.*

The proof of theorem 1 can be reduced to a new removability-type theorem to quasiconformal maps, Theorem 4: see the section 1.2 for the statement of this theorem, and, in the end of the that section, a short explanation of how this reduction is done in section 6, which can be read once the reader is familiar with the removability result and notation introduced in section 1.2.

Perhaps the assumptions about the maps $f_1$ and $f_2$ in Theorem 1 will sound quite artificial to the reader, specially the puzzle geometry control. But it seems that there are many combinatorial types of polynomial-like maps for which the critical puzzle pieces in the so called "principal nest" satisfies those assumptions. A well studied case is the Fibonacci combinatorics: this combinatorial type satisfies the assumptions for all even criticalities (for criticality two, see Lyubich[Lyu93]; for other criticalities, see Buff[Bu]). The result for criticalities larger than two will be very important in the proof of the hyperbolicity of the Fibonacci cycle, using Theorem 1 (see Theorem 2). For criticality two, recently Perez[P] proved that, given a hyperbolic quadratic polynomial, there exists a Cantor set of non renormalizable quadratic polynomials (but with zero Hausdorff dimension) in the boundary of the Mandelbrot set, so that the shape of the critical puzzle pieces in the principal nest converges (in the Hausdorff metric over compact sets) to the shape of the Julia set of this hyperbolic polynomial. This, together with the linear grow of the modulus in the principal nest, shows that induced maps in the principal nest of these quadratic polynomials satisfy the puzzle geometry control and the complex bounds in Theorem 1.

Two metric properties are now classic tools to prove rigidity: the "bounded geometry of the postcritical set" and the "linear decay of modulus": the first one was used to prove the rigidity of the infinitely renormalizable polynomials (with an unique critical point) with bounded combinatorics (see Sullivan[Sul92] and de Melo & van Strien[dMvS]). The second one was used in the proof of the real Fatou conjecture to quadratic polynomials (Lyubich[Lyu97] and Graczyk & Swiatek[GS97]).

It is quite remarkable that both behaviors can coexist with the puzzle geometry control: the Fibonacci combinatorics is the simplest example. For criticality two, Fibonacci maps have linear decay of modulus (Lyubich & Milnor[LM]); for larger even criticalities, the Fibonacci maps have a postcritical set with bounded geometry (see Bruin, Keller, Nowicki & van Strien[BKNS], for high criticalities, and Keller & Nowicki[KN]). So the "puzzle geometry control" property seems to be more "primitive" than the two classical properties above.

On the other hand, the linear decay of modulus holds for all real non renormalizable quadratic polynomials (Lyubich[Lyu94]), what is no longer true to the puzzle geometry control property. However, the linear decay of geometry lost its universality for larger criticalities, or even for multimodal maps with quadratic critical points (Swiatek & Vargas[SV]). An even worse fact is the possibility of existence of real non renormalizable combinatorics, for higher criticalities and/or many critical points, which does not have neither decay of modulus (at any rate) nor bounded geometry (or even essentially bounded geometry). Our hope is that the puzzle geometry control and complex bounds hold for some of these combinatorial types, whose rigidity is seemly impossible to prove using previously disponible methods (see section 17).

Apart this future development, we will shown how useful can be the Theorem 1, proving in the second part of this work that the Fibonacci renormalization operator has a hyperbolic circle. We are going to introduce the Fibonacci renormalization operator before to state the results.

After hyperbolic, infinitely renormalizable with bounded combinatorics and Collet-Eckmann polynomials with an unique critical point, the Fibonacci combinatorics is the combinatorial type with interesting properties which is better knowledge: We can, for instance, cite Lyubich & Milnor [LM]; Lyubich [Lyu93]; Bruin, Keller, Nowicki & van Strien [BKNS]; Keller & Nowicki [KN]; van Strien & Nowicki [vSN94]; Buff [Bu]. The main points are that

- The Fibonacci combinatorics is the simplest one which admits wild attractors,
- For critical order larger than two, the Fibonacci combinatorics have bounded geometry,
- The Fibonacci combinatorics is the simplest combinatorics which is infinitely renormalizable in the generalized sense.

Let $d$ be an even natural number larger than two ($d$ will be fixed for the rest of this work). If $\alpha$ is a complex number, and if $A$ is a subset of the complex plane, denote $\alpha A := \{\alpha \cdot x \colon x \in A\}$.

Let $f$ be a real analytic map defined in two intervals, a symmetric interval $I_0^1$ and $I_1^1$. Denote $I_0^0 := f(I_1^1)$. Assume that

- $f \colon I_0^1 \to I_0^0$ has zero as its unique critical point, which is a maximum. Moreover, the degree of zero is $d$ and $f(\partial I_0^1) \subset \partial I_0^0$.
- $f \colon I_1^1 \to I_0^0$ does not have critical points and $f(\partial I_1^1) \subset \partial I_0^0$.

We will say that $f$ is **Fibonacci renormalizable** if

- $f(0) \in I_1^1$,
- $f^2(0) \in I_0^1$,
- $f^3(0) \in I_0^1$,
- The connected component of the domain of $f^2$ that contains the critical point also contains a fixed point, denoted $\beta_f$, so that $Df^2(\beta_f) < 0$.

Denote by $I_0^2$ the connected component of the domain of the first return map to $I_0^1$ which contains 0 and by $I_1^2$ the connected component of the domain of the first return map to $I_0^1$ which contains $f^2(0)$. Note that $I_0^2 \cap I_1^2 = \phi$. The **Fibonacci renormalization** of $f$, denoted $\mathcal{R}f$, is defined as the map

$$\mathcal{R}f \colon \frac{1}{\beta_f}I_0^2 \cup \frac{1}{\beta_f}I_1^2 \to \frac{1}{\beta_f}I_0^1,$$

where $\mathcal{R}f(x)$ is equal to $1/\beta_f f^2(\beta_f \cdot x)$ on $1/\beta_f I_0^2$, and $1/\beta_f f(\beta_f \cdot x)$ on $1/\beta_f I_1^2$. If $\mathcal{R}^n f$ is defined for every $n \geq 0$, we say that $f$ is **infinitely renormalizable** in the Fibonacci sense.

In van Strien & Nowicki[vSN94], it was proved that this operator have an unique periodic orbit of period two, the **Fibonacci cycle**. Indeed, if $\tilde{f}_1 \colon I \cup J \to T$ is one of the maps in the periodic orbit, then the map $\tilde{f}_2$ defined as $\tilde{f}_1$ on $I$ and $-f_1$ on $J$ is exactly $\mathcal{R}f$ (see Buff[Bu]). Making a slight modification of the Fibonacci renormalization operator, we can consider these maps as fixed points (see section 16). The following remark is very important: Since these maps are almost the same, a distinction between the two maps will be not important in some parts of this work, so sometimes we will denote both of them by $\tilde{f}$. For instance $|g - \tilde{f}|$ denotes $\min\{|g - \tilde{f}_1|, |g - \tilde{f}_2|\}$.

Indeed, as proved by van Strien & Nowicki[vSN94], $\tilde{f}$ has a generalized polynomial-like extension $\tilde{f} \colon \overline{U}_0 \cup \overline{U}_1 \to \overline{V}$, where $0 \in U_0$. Selecting $U_0$ and $U_1$ in an appropriated way, we will define an operator $\mathcal{R} \colon B_{\epsilon_0}(\tilde{f}) \to \mathcal{B}_{nor}(U_0 \cup U_1)$ (this is an affine subspace of complex analytic functions, with an appropriated normalization, with continuous extension to $\overline{U}_0 \cup \overline{U}_1$), where $B_{\epsilon_0}(\tilde{f}) \subset \mathcal{B}_{nor}(U_0 \cup U_1)$ is the $\epsilon_0$-ball around $\tilde{f}$, which coincides with the Fibonacci renormalization operator above defined in the real trace of the real functions in $B_{\epsilon_0}(\tilde{f})$, so that:

**Theorem 2** (Hyperbolicity). *The Fibonacci renormalization operator $\mathcal{R}$, defined in a small neighborhood of $\tilde{f}$ in $\mathcal{B}_{nor}(U_1 \cup U_2)$, is a compact operator with a hyperbolic cycle $\{\tilde{f}_1, \tilde{f}_2\}$ with codimension one stable manifold $W^s(\tilde{f})$. Furthermore, if $f$ is real analytic map that is infinitely renormalizable in the Fibonacci sense, then there exists $N = N(f)$ so that $\mathcal{R}^N f \in W^s(\tilde{f})$.*

One of the ingredients to prove Theorem 2 is the control of the shapes of the critical puzzle pieces in the principal nest of the maps in the Fibonacci cycle: Using the Small Orbits Theorem[Lyu99] and the contraction in the hybrid class, the proof of the hyperbolicity is reduced to a rigidity problem: given a map $f$ whose all Fibonacci renormalizations are very close of the Fibonacci cycle, we need to construct a quasiconformal conjugacy between this map and one of the maps in the Fibonacci cycle. The first observation is that it is possible to control the geometry of the correspondent puzzle pieces of $f$ (see the end of section 9). The next step is to prove that there are a topological conjugacy between $f$ and one of the maps in the Fibonacci cycle which is quasiconformal outside the Julia of $f$. Finally Theorem 1 is used to conclude that the conjugacy is, indeed, quasiconformal everywhere.

The proof of the exponential contraction on the hybrid classes of the maps in the Fibonacci cycle uses a new argument (see section 14). Once we have the non existence of invariant line fields in the Fibonacci tower and the (topological) convergence of the Fibonacci renormalization operator in the hybrid classes (both of them are van Strien & Nowicki results [vSN94]), the proof is infinitesimal: we prove that the derivative of the Fibonacci operator is a contraction in the 'tangent spaces' of the hybrid classes. The proof is more simple and straightforward than the previous arguments in McMullen[McM] and Lyubich[Lyu99], but indeed it is not necessary any breakthrough, so the reader must consider it as a natural development of McMullen and Lyubich proofs. However, the method seems to be so general as the previous ones: it can be applied in the classic renormalization horseshoe [Lyu99], for instance.

The hyperbolicity of the Fibonacci renormalization operator implies, after a short argument in section 16, the following universality result:

**Theorem 3** (Universality). *Let $d$ be a even natural number larger than two. Then there exists $\gamma = \gamma_d > 1$ so that, for an open and dense set of families $f_\lambda \in \mathcal{B}_{nor}(U)$ of even generalized polynomial-like maps whose domain contains two connected components, which is real for real parameters, the following holds: Let $f_{\lambda_\infty}$ be quasi-conformally conjugated with either $\tilde{f}_1$ or $\tilde{f}_2$. Then there exists $\delta > 0$, $N \in \mathbb{N}$, and a sequence $\lambda_n$ so that $\lambda_n \to \lambda_\infty$ and*

- *For $n \geq N$, $\lambda_n$ is the unique parameter so that $|\lambda_n - \lambda_\infty| \leq \delta$ and whose map $f_{\lambda_n}$ has a superattractor with period $S_n$ (the Fibonacci sequence) and whose closest returns of the critical point are exactly $S_i$, $i < n$.*
- *We have*
$$\frac{|\lambda_n - \lambda_{n+1}|}{|\lambda_{n+1} - \lambda_{n+2}|} \to \gamma.$$

*Furthermore, the same statement holds to the family $x^d + c$ (perturbations are not necessary, in this case).*

1.2. **Removability result.** We will state a new abstract removability-type result: this result is the key to Theorem 1. Before to state it, we need to introduce some trivial notation. Assume in definitions that all functions are measurable.

In this section, we will reserve the word **family** to denote families of domains: a domain is an open and connected set in $\mathbb{C}$. Moreover, we will assume that all these domains have boundaries with zero 2-dimensional Lebesgue measure. Sometimes we will deal with families of **centered** domains: this is a family of pairs $(A, x)$, where $x \in A$. We will say that $x$ is a center of $A$. Note that $A$ can have many centers and $x$ can be a center of many domains. Recall that $m(A)$ denotes the 2-dimensional Lebesgue measure of $A$.

**Definition 1.2.** *A family $\mathcal{F}$ covers a set $X \subset \mathbb{C}$ in **arbitrary small scales** if for any $x \in X$ and $\epsilon > 0$ there exists $A \in \mathcal{F}$ such that $x \in A$ and $\text{diam } A \leq \epsilon$. Moreover $\partial A \cap X = \phi$ for $A \in \mathcal{F}$.*

**Definition 1.3.** *A centered family $\mathcal{F}$ covers a set $X \subset \mathbb{C}$ in **all scales** if for any $x \in X$ there exists $C = C(x) > 1$ so that, for any $\epsilon \leq \epsilon_0$, with $\epsilon_0 = \epsilon_0(x)$, $x$ is the center of one domain $A \in \mathcal{F}$ satisfying:*
$$\frac{1}{C}\epsilon \leq \text{diam } A \leq C\epsilon.$$

*Moreover $X \cap \partial A = \phi$ for $A \in \mathcal{F}$. If $C$ and $\epsilon_0$ does not depend on $x$, we say that $\mathcal{F}$ covers a set $X$ in all scales in an $(C, \epsilon_0)$-uniform way.*

**Definition 1.4.** *We say that a family $\mathcal{F}$ has **uniform bounded geometry** if there exists $\lambda < 1$ so that for every $A \in \mathcal{F}$ there exists $x \in A$ so that*
$$B(x, \lambda \cdot \text{diam } A) \subset A.$$

**Definition 1.5.** *We say that a centered family $\mathcal{F}$ has **pointwise bounded geometry** if for any $x \in \mathbb{C}$ there exists $C = C(x) < 1$ so that if $(A, x) \in \mathcal{F}$ then $B(x, C \cdot diam\ A) \subset A$. The family $\mathcal{F}$ has $C$-**uniform bounded geometry** if $C$ does not depend on $x$.*

**Definition 1.6.** *We say that a family of domains $\mathcal{F}$ is a **Markov family** if*
- *If $A_1$ and $A_2$ belongs to $\mathcal{F}$, then one of the following statements holds:*
  - *$A_1 \cap A_2 = \phi$,*
  - *$A_1 \subset A_2$,*
  - *$A_2 \subset A_1$.*
- *Suppose that $\mathcal{F}' \subset \mathcal{F}$ is a subfamily with the following properties:*
  - *The inclusion is a total order in $\mathcal{F}'$: if $A_1, A_2 \in \mathcal{F}'$, then either $A_1 \subset A_2$ or $A_2 \subset A_1$.*
  - *There exists $A_0 \in \mathcal{F}$ such that $A_0 \subset A$ for any $A \in \mathcal{F}'$.*
  
  *then $\mathcal{F}'$ is finite.*

Markov families arise in a very natural way in complex dynamics: the family of the puzzle pieces of any level is a Markov family. Note that any subfamily of a Markov family is a Markov family itself. The main property of Markov families is the following: if a Markov family $\mathcal{F}$ covers a set $X$, then there exists a subfamily $\mathcal{F}' \subset \mathcal{F}$ which covers each point of $X$ once and almost each point of $\mathbb{C}$ at most once. The reader can fill the proof of this property.

Let $h$ be a homeomorphism in the complex plane and let $\mathcal{F}$ be a (centered) family. Define
$$h(\mathcal{F}) = \{h(A) \colon\ A \in \mathcal{F}\},$$
if $\mathcal{F}$ is a family, or
$$h(\mathcal{F}) = \{(h(A), h(x)) \colon\ (A, x) \in \mathcal{F}\},$$
if $\mathcal{F}$ is a centered family.

**Theorem 4** (Removability Theorem)**.** *Let $h \colon \mathbb{C} \to \mathbb{C}$ be a homeomorphism and let $X \subset \mathbb{C}$ be a compact set. Suppose that $h$ is quasiconformal in $\mathbb{C} - X$. Assume that there are a decomposition $X = X_1 \cup X_2$, a family $\mathcal{F}_1$ and a centered family $\mathcal{F}_2$ such that, with respect to $X_1$ and $\mathcal{F}_1$:*
- *$\mathcal{F}_1$ is a Markov family,*
- *$\mathcal{F}_1$ covers $X_1$ in arbitrary small scales,*
- *$\mathcal{F}_1$ and $h(\mathcal{F}_1)$ have uniform bounded geometry.*

*and, with respect to $X_2$ and $\mathcal{F}_2$:*
- *$X_2$ has zero 2-dimensional Lebesgue measure,*
- *$\mathcal{F}_2$ covers $X_2$ in all scales,*
- *$\mathcal{F}_2$ and $h(\mathcal{F}_2)$ have pointwise bounded geometry.*

*Then $h$ is quasiconformal in $\mathbb{C}$.*

**Remark 1.4.** *We can add another exceptional set $X_3$ so that $X_3$ has $\sigma$-finite 1-dimensional Hausdorff measure and anything more about $X_3$ is assumed.*

**Remark 1.5.** *The same proof works replacing $\mathbb{C}$ by $\mathbb{R}^n$, $n \geq 2$, and quasiconformal by quasisymmetric, assuming that $X_2$ has zero $n$-dimensional Lebesgue measure and $X_3$ has $\sigma$-finite $(n\text{-}1)$-dimensional Hausdorff measure.*

Theorem 4 is an almost straightforward generalization of the main result in Kallunki & Koskela[KK00] and so its proof will not be hard for the reader familiar with the methods there: a special property about certain coverings by closed discs, the so called Besicovitch Theorem, is replaced in some parts of the proof by our Markov property. Replace discs by (centered) domains with bounded geometry in some points of the proof will not cause much trouble. Anyway, for completeness we present an almost self contained proof in section 5. Apart the above definitions about families, section 5 is independent of the rest of this work.

In the application of Theorem 4 to prove Theorem 1 (see section 6), $X$ will be the filled-in Julia set (with empty interior) of a generalized polynomial-like map, $X_1$ will be the points in the Julia set whose forward orbit accumulates in the critical point and $X_2$ will be the set of points whose forward orbit is far of the critical point. The family $\mathcal{F}_1$ will be defined as the family of univalent pullbacks of critical puzzle pieces with bounded geometry. The complex bounds and the Koebe distortion lemma will provide the uniform bounded geometry. For a given point in $X_2$, the domains in $\mathcal{F}_2$ which cover this point will be univalent pullbacks with bounded distortion ( but whose control depends on the point) of a finite set of puzzle pieces (which depends on the point, too) which do not contain the critical point. This will give us the pointwise bounded geometry.

## 2. Outline of the work

In the first part of this work we prove Theorem 4 and Theorem 1. We prove Theorem 4 in section 5 and we reduce the Theorem 1 from Theorem 4 in section 6. A simple application of Theorem 1 is given in section 7: we give a new proof of the rigidity of the Fibonacci combinatorics in the real quadratic family. In the second part we prove the hyperbolicity of the Fibonacci cycle and the universality result. In section 8 we study the maps in the Fibonacci cycle, using results in Buff[Bu]. A complex analytic version of the Fibonacci renormalization operator is defined in section 9. In section 10, we prove that the geometry of certain critical puzzle pieces of maps whose renormalizations are very close of the Fibonacci cycle is under control: this will be critical to obtain the hyperbolicity of the Fibonacci cycle using Theorem 1. In section 11 we prove that, if the orbit of a map by the Fibonacci renormalization operator is always very close to the Fibonacci cycle, then there exists a hybrib conjugacy between that map and one of the maps in the cycle. In section 12 we give a introduction to the theory of infinitesimal perturbations of maps developed in Lyubich[Lyu99] and Avila, Lyubich & de Melo[ALdM], with minor modifications, and in section 13 we introduce some basic facts about induced transformations: the Fibonacci renormalization operator is an example of induced transformation. In section 14, the exponential contraction of the Fibonacci renormalization operator in the hybrid class of the maps in the Fibonacci cycle is proved, using the non existence of invariant line fields in the Fibonacci tower[vSN94] and a new simple argument. Finally, in section 15 we prove that the Fibonacci cycle is hyperbolic (Theorem 2). In section 16 we derive the universality result (Theorem 3). Future developments and open questions are described in section 17.

## 3. Motivation

The initial motivation to Theorem 1 was the study of the Fibonacci renormalization (to high criticalities): by Small Orbits Theorem [Lyu99] and the exponential contraction on the hybrid classes, to prove that the Fibonacci cycle is hyperbolic, it is sufficient to prove that maps whose all Fibonacci renormalizations are very close to the Fibonacci cycle are indeed quasiconformally conjugated to one of the maps in the Fibonacci cycle. So the problem was reduced to a rigidity result. The first step is to find a conjugacy which is quasiconformal outside the filled-in Julia sets. This is done in section 11.

The next step is to prove that the conjugacy is indeed quasiconformal in the Julia set. Unfortunately, all the previous methods to prove this kind of result to non-renormalizable (in the classical sense) and non necessarily real maps are based on the "modulus divergence property", "linear decay of moduli property" or that the dynamics has some "expansivity" (Collet-Eckmann, somability or Rivera-Letelier's decay of geometry). Any of these properties does not hold for the Fibonacci maps with criticalities larger than two. Furthermore, in all these proofs the quasiconformal removability, dynamical removability or at least zero Lebesgue measure of the Julia set was derived. Similar questions for infinitely classic-renormalizable maps are major problems in complex dynamics, but to prove the hyperbolicity of the classic renormalization operator, it is not necessary to solve them: there are not reason to believe that the situation is distinct for the Fibonacci renormalization.

So a new ingredient seems to be necessary: recent works on the Fibonacci map (see van Strien & Nowicki [vSN94] and Buff [Bu]) convinced ourselves that the control of the puzzle's geometry is such ingredient: the boundaries's shapes of the critical puzzle pieces in the principal nest of a real Fibonacci generalized polynomial-like maps converges to the shape of Julia set of a hyperbolic polynomial-like map (proved by Buff). We observed that this property is robust, in the sense that it is very easy to transfer a similar property to a map whose renormalizations are very close of the maps in the Fibonacci cycle: the correspondent boundaries's puzzles will not (a priori) necessarily converges to some shape, but its geometry will be under control (see section 9). But how to use this property to prove the wished rigidity result? The naive argument is the following: we have, by assumption, that some conjugacy $h$ does not distort certain critical puzzles pieces. Using the dynamics, we transport this control of puzzle distortion to any point whose forward orbit accumulate in the critical point (here it is necessary some "complex bounds"). We can deal with the points whose orbit does not accumulate the critical point in another way (see below). So we conclude that the conjugacy does not distort a lot of puzzle pieces, around a lot of points in the Julia set (this is done in section 6): we can expect that this conjugacy is very regular.

The recent proof by Przytycki and Rohde [PR99] that the Julia set of CE holomorphic repellers are dynamically removable awake our attention to a result by Heinonen & Koskela [HK95], which claims that, to prove that a homeomorphism is quasiconformal, it is sufficient to find, for each point in the domain, a sequence of circles centered on this point, whose diameter goes to zero, so that the geometry of the image of these circles by the homeomorphism is under control. Since we will not use circles, but puzzle pieces (which hardly are circles), a new kind to "removability" result is necessary: certain coverings of sets by discs have very special properties, as the Besicovitch theorem. But the advantage of the puzzle pieces is that they satisfies a "Markov property" (see definition 1.6), which substitute the Besicovich theorem, and they are dynamically defined, what is very unusual for circles!

Unfortunately, we cannot transport the information about the geometry the critical puzzle pieces in the principal nest and their images to all points in the Julia set, since there are points whose orbits do not accumulate the critical point. The solution to this problem is that the set of these problematic points is small (zero Lebesgue measure), so we can be (in some sense) less exigent with them: we need only a pointwise control of distortion, which is obtained using the expansion of the dynamics far of the critical point. Here a argument developed by Kallunki & Koskela [KK00] will be an useful tool. So a new removability-type theorem arise (Theorem 4).

Note that the hyperbolicity of the Fibonacci renormalization operator does not exhaust the applications of Theorem 1. See section 17.

## 4. Preliminaries

4.1. **Quasiconformal theory.** We describe here some facts about quasiconformal maps and vectors. However we will assume that the reader is familiar with the general properties of quasiconformal maps (Lehto & Virtanen[LV] is a good reference). Let $f$ be a function whose distributional derivatives belongs to $L^2_{loc}$. Define

$$\partial f = f_z = \frac{1}{2}(f_x - if_y).$$

$$\overline{\partial} f = f_{\overline{z}} = \frac{1}{2}(f_x + if_y).$$

By the Weil lemma, if $f_{\overline{z}}$ is zero then $f$ is an analytic function.

A vector field $v$ in the Riemann sphere is a **K-quasiconformal vector field** if $|v_{\overline{z}}|_{\overline{\mathbb{C}}} < K < \infty$. Here $|\cdot|_{\overline{\mathbb{C}}}$ denotes the spherical metric. The following property will be very useful:

**Proposition 4.1** ([McM]). *Let $z_0$, $z_1$ and $z_2$ be distinct points in the Riemann sphere. and $K > 0$. Then the set of $K$-quasiconformal vector fields in the Riemann sphere which vanishes at these points is a compact set in the space of continuous vector fields in the Riemann sphere with the sup norm induced by the spherical metric.*

**Proposition 4.2** (Lemma 10 in [AB]). *If $f$ and $g$ are functions whose distributional derivatives belong to $L^2_{loc}$ then $f \circ g$ also have derivatives in $L^2_{loc}$ and*

$$(f \circ g)_z = (f_z \circ g)g_z + (f_{\overline{z}} \circ g)\overline{g}_z.$$

$$(f \circ g)_{\overline{z}} = (f_z \circ g)g_{\overline{z}} + (f_{\overline{z}} \circ g)\overline{g}_{\overline{z}}.$$

The **Beltrami field** associated to a homeomorphism $f$ with derivatives in $L^2_{loc}$ is

$$\mu_f = \frac{\partial_{\overline{z}} f}{\partial_z f}.$$

If $|\mu_f|_\infty < \infty$ then we call $f$ a **quasiconformal** homeomorphism. If $f$ and $g$ are quasiconformal homeomorphism then $f \circ g$ is a quasiconformal homeomorphism with the following Beltrami field (see, for instance,[LV]):

$$\mu_{f \circ g^{-1}}(g(z)) = \frac{\mu_f(z) - \mu_g(z)}{1 - \mu_f(z)\overline{\mu_g(z)}}\left(\frac{g_z(z)}{|g_z(z)|}\right)^2.$$

Let $f$ be an analytic function. The **pullback** of a Beltrami field $\mu$ by $f$ is defined by

$$f_*\mu(z) = \mu(f(z))(\frac{f_z}{|f_z|})^{-2}.$$

The Beltrami field $\mu$ is **invariant** by $f$ if $f_*\mu = \mu$ (in the intersections of the domains of definition of $\mu$ and $f_*\mu$). A quasiconformal vector field $v$ is a **deformation** of $f$ if $\overline{\partial}v$ is invariant by $f$.

4.2. **Banach spaces.** We will denote domains (pre-compact open sets in the complex plane) with symbols like $U, V, \tilde{V}, \hat{U}, \ldots$. If we denote a domain as, for instance, $U = U_1 \cup \cdots \cup U_n$, the reader must assume that $U_i$ are the (simply) connected components of $U$. If $\alpha \in \mathbb{C}$, then $\alpha U := \{\alpha \cdot x \colon x \in U\}$. The $\delta$-neighborhood of a set $X$, denoted $\delta$-X, is the set $\{x \colon dist(x, X) \leq \delta \; diam \; X\}$. Let $U$ be a domain so that $\overline{U}$ does not disconnect the plane. Let $d > 1$. In the Banach space $\mathcal{B}(U)$ of holomorphic functions on $U$ with continuous extension to $\overline{U}$, provided with the sup norm, consider the closed affine space $\mathcal{B}_{nor}(U)$ of functions $f$ so that

- the domain $U$ contains 1 and 0 and furthermore $f(1) = 1$ and $f^{(i)}(0) = 0$, for $1 \leq i < d$.

We will call $\mathcal{B}_{nor}(U)$ a **Banach slice**. The tangent space of $\mathcal{B}_{nor}$, denoted by $T\mathcal{B}_{nor}(U)$, is the subspace of the vectors $v \in \mathcal{B}(U)$ so that

- $v(1) = 0$,
- $v^{(i)}(0) = 0$, for $1 \leq i < d$.

Let $f \colon U_0 \cup \cdots \cup U_n \to V$, where $U_i$ are simply connected, be a generalized polynomial-like map and let $\hat{U} = \tilde{U}_0 \cup \cdots \cup \tilde{U}_n$, where $\tilde{U}_i$ are simply connected, be a domain so that $\tilde{U}_i \Subset U_i$. Then we say that the Banach space $\mathcal{B}_{nor}(\tilde{U})$ is **compatible** with $f$. By an Levin & van Strien result (proposition 2.2 in [LvS]), if $g \colon \hat{U}_1 \cup \cdots \cup \hat{U}_n \to \hat{V}$ is a generalized polynomial-like restriction of $f$, with $\hat{U}_i \subset \tilde{U}_i$, then $K(g) = K(f)$. We will say that $g$ is a **representation** of $f$ subordinated to $\tilde{U}$. So all these representations have the same filled-in Julia set and, if the critical points does not escape for $f$, the same postcritical set.

4.3. **The principal nest.** Let $f \colon \cup_i U_i \to V$ be a generalized polynomial-like map whose critical point does not escape. Denote $V_0^0 := V$, $V_0^1 := U_0$ and $V_1^1 := U_k$, where the critical point is in $U_0$ and the critical value in $U_k$. Denote $f_0 := f$. $f$ restricted to $V_0^1$ will be denoted $f_0^0$ and $f$ restricted to $V_1^1$ by $f_0^1$. denote $\mathcal{E}_0 := \phi$ and $\mathcal{D}_0 := \{V_0^1, V_1^1\}$. Suppose, by induction, that we have defined $V_0^n$, with $0 \in V_0^n$ and the following property:

$$f^i(\partial V_0^n) \cap V_0^n = \phi,$$

for every $i > 0$ so that $f^i$ is defined on $V_0^n$. Denote by $\mathcal{E}_n$ the family of (simply) connected components of

$$\{x \in \cup_i U_i \colon \text{ there exists } m \geq 0 \text{ so that } f^m(x) \in V_0^n\}.$$

The **first entry map** to $V_0^n$ of a point $x \in \cup \mathcal{E}_n$, denoted by $\pi_n(x)$, is defined as $f^m(x)$, where $m \geq 0$ is minimal so that $f^m(x) \in V_0^n$. The **first return map** to $V_0^n$, denoted by $f_n$, is defined as $\pi_n \circ f$ in $\cup \mathcal{D}_n$, where $\mathcal{D}_n$ is the family of (simply) connected components of $f^{-1}E$, where $E \in \mathcal{E}_n$. Note that $\mathcal{E}_n - \{V_0^n\} \subset \mathcal{D}_n$. Denote by $V_0^{n+1} \in \mathcal{D}_n$ the domain which contains the critical point and by $V_1^{n+1} \in \mathcal{D}_n$ the domain which contains $f_n(0)$. Note that $V_0^n$ is a nested sequence of critical puzzle pieces, called **principal nest**.

## Part 1. Puzzle Geometry and Rigidity

### 5. Removability theorem

We are going to prove Theorem 4. The following proof is almost self-contained. We follow closely the ideas of S. Kallunki & P. Koskela([KK00]) and J. Heinonen & P. Koskela([HK95]). Denote $B(z_0, r) = \{z\colon |z - z_0| \leq r\}$, $S(z_0, r) = \{z\colon |z - z_0| = r\}$ and define $\pi_{z_0}\colon \mathbb{C} - \{0\} \to S(z_0, 1)$ by $\pi_{z_0}(z) = \frac{z - z_0}{|z - z_0|} + z_0$. Denote by $1_A$ the characteristic function of the set $A$. $m$ is the 2-dimensional Hausdorff measure and $m_1$ is the 1-dimensional Hausdorff measure.

*Proof of Theorem 4.* Let $x_0$ be an arbitrary point and $r_0 > 0$. Define

$$l_h(x_0, r_0) := \inf\{|h(x) - h(x_0)|\colon |x - x_0| \geq r_0\}.$$

$$L_h(x_0, r_0) := \sup\{|h(x) - h(x_0)|\colon |x - x_0| \leq r_0\}.$$

We need to prove that

$$\frac{L_h(x_0, r_0)}{l_h(x_0, r_0)} \leq H.$$

for some universal $H$. We will fix $x_0$ and $r_0$. Denote $L = L_h(x_0, r_0)$ and $l = L_h(x_0, r_0)$. Define $A := h^{-1}(\{y\colon l \leq |y - h(x_0)| \leq L\}) - X_2$. Since $h$ is quasiconformal outside $X$, there exists $\lambda < 1$ so that, for every $x$ outside $X$, we can find $r_0(x) > 0$ so that

$$B(h(x), \lambda \ diam \ h(B(x, r))) \subset h(B(x, r)),$$

for $r < r_0(x)$. Let $j_0$ be minimal such that $L2^{-j_0} < l$. We will define a family of closed discs $\mathcal{B}_j$ and a Markov family $\mathcal{M}_j$ in the following way:

- A disc $B(x, R)$ belongs to $\mathcal{B}_j$, with $j \leq j_0 + 1$, if
  - $x \in A - X$,
  - $r < r_0(x)$ and
  - $h(B(x, r))$ is contained in $R_j := \{y\colon L2^{-(j+1)} < |y - h(x_0)| < L(2^{-j} + 2^{-(j_0+1)})\}$ and $diam\ h(B(x, r)) \leq L2^{-100j_0}$,
- A piece $M \in \mathcal{F}_1$ belongs to $\mathcal{M}_j$, with $j \leq j_0 + 1$, if $h(M)$ is contained in $R_j$ and $diam\ h(M) \leq L2^{-100j_0}$.

Note that $\cup_j \mathcal{B}_j$ covers $A - X$ and $\cup_j \mathcal{M}_j$ covers $A \cap X_1$. Applying the Besicovitch's covering theorem to $\cup_j \mathcal{B}_j$, we can

- Replace $\mathcal{B}_j$ by an enumerable subfamily, denoted by the same symbol $\mathcal{B}_j$, such that

$$1_{A-X} \leq \sum_{j \leq j_0} \sum_{B \in \mathcal{B}_j} 1_B \leq C.$$

- Moreover, each family $\mathcal{B}_j$ can be decomposed in $Q$ subfamilies $\mathcal{B}_j^i$, $i = 1, \ldots, Q$, where $Q$ does not depends on anything, so that

$if\ B_1 \in \mathcal{B}_{j_1}^i\ and\ B_2 \in \mathcal{B}_{j_2}^i\ then\ either\ B_1 \cap B_2 = \phi, or\ j_1 = j_2\ and\ B_1 = B_2.$

Similarly, since $\mathcal{F}_1$ is a Markov family, we can replace $\mathcal{M}_j$ by a subfamily, if necessary, so that
$$\sum_{j \leq j_0} \sum_{M \in \mathcal{M}_j} 1_M \leq 1$$
almost everywhere in $\mathbb{C}$, and
$$\sum_{j \leq j_0} \sum_{M \in \mathcal{M}_j} 1_M = 1$$
in $A \cap X_1$. Note that there exists an universal constant $C > 1$ so that, for any $M \in \cup_j \mathcal{M}_j$, there exists a disc $B^M$ which satisfies
- $diam\ B^M \sim diam\ M$,
- $diam\ C\ B^M > 2\ diam\ M$,
- $B^M \subset M \subset C\ B^M$.

Define
$$\rho := \left(\log \frac{L}{l}\right)^{-1} \sum_{j \leq j_0} \sum_{B \in \mathcal{B}_j} \frac{diam\ h(B)}{dist(h(B), h(x_0))} \frac{1}{diam\ B} 1_{2B}$$
$$+ \left(\log \frac{L}{l}\right)^{-1} \sum_{j \leq j_0} \sum_{M \in \mathcal{M}_j} \frac{diam\ h(M)}{dist(h(M), h(x_0))} \frac{1}{diam\ M} 1_{CB^M}.$$

We claim that
$$\int \rho^2\ dm \leq \left(\log \frac{L}{l}\right)^{-1}.$$

Denote
$$a_B := \frac{diam\ h(B)}{dist(h(B), h(x_0))} \frac{1}{diam\ B}.$$
$$a_M := \frac{diam\ h(M)}{dist(h(M), h(x_0))} \frac{1}{diam\ M}.$$

By lemma 4.2 in [B],
$$\int \rho^2\ dm = \left(\log \frac{L}{l}\right)^{-2} \left(\int \left(\sum_{j \leq j_0} \sum_{B \in \mathcal{B}_j} a_B 1_{2B} + \sum_{j \leq j_0} \sum_{M \in \mathcal{M}_j} a_M 1_{CB^M}\right)^2 dm\right)$$
$$\leq C \left(\log \frac{L}{l}\right)^{-2} \left(\int \left(\sum_{j \leq j_0} \sum_{B \in \mathcal{B}_j} a_B 1_B\right)^2 dm + \int \left(\sum_{j \leq j_0} \sum_{M \in \mathcal{M}_j} a_M 1_{B^M}\right)^2 dm\right)$$
$$\leq C \left(\log \frac{L}{l}\right)^{-2} \left(\sum_{j \leq j_0} \sum_{B \in \mathcal{B}_j} a_B^2 m(B) + \sum_{j \leq j_0} \sum_{M \in \mathcal{M}_j} a_M^2 m(M)\right).$$

by assumption
$$\frac{m(h(B))}{(diam\ h(B))^2},\ \frac{m(h(M))}{(diam\ h(M))^2},\ and\ \frac{m(M)}{(diam\ M)^2}$$
are at a definitive distance of zero and infinity. We obtain
$$\int \rho^2\ dm$$
$$\leq C \left(\log \frac{L}{l}\right)^{-2} \left(\sum_{j \leq j_0} \sum_{B \in \mathcal{B}_j} \frac{m(h(B))}{dist(h(B), h(x_0))^2} + \sum_{j \leq j_0} \sum_{M \in \mathcal{M}_j} \frac{m(h(M))}{dist(h(M), h(x_0))^2}\right)$$
$$\leq C \left(\log \frac{L}{l}\right)^{-2} \left(\sum_{j \leq j_0} (2^{-j}L)^{-2} \sum_{B \in \mathcal{B}_j} m(h(B)) + \sum_{j \leq j_0} (2^{-j}L)^{-2} \sum_{M \in \mathcal{M}_j} m(h(M))\right)$$

$$\leq C \left(log \frac{L}{l}\right)^{-1}.$$

The last inequality is consequence of the following:
- $h(\mathcal{B}_j)$ covers each point in $R_j$ at most an universal number of times,
- $h(\mathcal{M}_j)$ covers almost each point in $R_j$ at most once,
- The measure of $R_j$ is commensurable with $(2^{-j}L)^2$,
- Assume that $j_0 > 10$: otherwise there are anything to prove. Then $j_0$ is smaller than $C \ log\frac{L}{l}$.

Without loss of generality (see [KK00]), we can assume that there are subcontinuas $F_1 \subset h^{-1}\{y\colon |y - h(x_0)| = l\}$ and $F_2 \subset h^{-1}\{y\colon |y - h(x_0)| = L\}$ such that $dist(F_1, F_2) = 2diam\ F_1 = 2diam\ F_2 = \frac{1}{2}$. Select $x_i \in F_i$ so that $|x_1 - x_2| = dist(F_1, F_2)$ and denote by $S$ the unique segment of length one whose middle point is the intersection point of $S$ and the segment between $x_1$ and $x_2$ and, moreover, $S$ cuts this segment exactly in the middle, in a perpendicular way. Using Lemma 2.1 in [KK00], to get a lower bound to $\int \rho^2\ dm$, it is sufficient to prove that, for fixed $x \in F_1$ and $y \in F_2$,

$$\int_{\gamma(x,y,z)} \rho\ dm \geq \delta > 0$$

for almost every $z \in S$. Here $\gamma(x, y, z)$ is the union of the segments between $x$ and $z$ and $y$ and $z$. Indeed,

$$\int_{\gamma(x,y,z)} \rho\ dm \geq C\left(\log\frac{L}{l}\right)^{-1} \sum_{j \leq j_0} \left(2^{-j}L\right)^{-1} \sum_{B \in \mathcal{B}_j,\ B \cap \gamma(x,y,z) \neq \phi} diam\ h(B)$$

$$+ C\left(\log\frac{L}{l}\right)^{-1} \sum_{j \leq j_0} \left(2^{-j}L\right)^{-1} \sum_{M \in \mathcal{M}_j,\ M \cap \gamma(x,y,z) \neq \phi} diam\ h(M).$$

Assume for a moment that $m_1(h(\gamma(x,y,z)) \cap h(X_2)) = 0$. Then

$$\sum_{B \in \mathcal{B}_j,\ B \cap \gamma(x,y,z) \neq \phi} diam\ h(B) + \sum_{M \in \mathcal{M}_j,\ M \cap \gamma(x,y,z) \neq \phi} diam\ h(M) \geq C\ 2^{-j}L.$$

Since $j_0 \sim \log \frac{L}{l}$, we obtained the lower estimation to $\int_{\gamma(x,y,z)} \rho\ dm$.

So we need to prove that for $x$ and $y$ fixed, $m_1(h(\gamma(x,y,z)) \cap h(X_2)) = 0$ for almost every $z \in S$. Indeed, fix $z_0 \in \mathbb{C}$ and denote $E_\sigma = \{z_0 + \lambda\sigma\colon 1 \leq \lambda \leq 2\} \cap X_2$. We claim that $m_1(h(E_\sigma)) = 0$, for almost every $\sigma \in \mathbb{S}^1$. Since $m(X_2) = 0$, we have that, for almost every $\sigma \in \mathbb{S}^1$, $m_1(E_\sigma) = 0$. We can decompose $E_\sigma$ in a countable number of parts ( denote each part by $E_\sigma$, again), so that for each part there exist constants $C_1 > 1$, $C_2, C_3 < 1$ and $r_0 > 0$ so that for any $x \in E_\sigma$ and $r < r_0$, there exists $(A, x) \in \mathcal{F}_2$ satisfying

- $\frac{1}{C_1}\ r \leq diam\ A \leq C_1\ r$,
- $B(x, C_2\ diam\ A) \subset A$,
- $B(h(x), C_3\ diam\ h(A)) \subset h(A)$.

By a simple argument as in the proof of lemma 2.4 in [KK00], it is sufficient to prove the following: there exists $C = C(\sigma, C_1, C_2, C_3) > 0$ so that if $F$ is a compact set in $z_0 + \sigma\mathbb{R}$, then $m_1(h(F \cap E_\sigma)) \leq C(m_1(F))^{1/2}$. To prove this, find intervals of disjoint interior $I_i \subset z_0 + \sigma\mathbb{R}$, $i = 1, \ldots, p$, with length $r \leq C_2/C_1 r_0$ so that $F \subset \cup_i I_i$ and $pr \leq 2\ m_1(F)$. For each $i$ so that $E_\sigma \cap I_i \neq \phi$, select $x_i \in I_i \cap E$. Then there exists $(A_i, x_i) \in \mathcal{F}_2$ so that

- $I_i \subset A_i$,
- $\frac{1}{C} r \leq diam A_i \leq Cr$.

So $h(F) \subset \cup_i h(A_i)$. Note that the family $\{h(A_i)\}_i$ covers each point in $\mathbb{C}$ at most $N(C_1, C_2)$ times. Then

$$\Big(\sum_i diam\ h(A_i)\Big)^2 \leq Cp \sum_i (diam\ h(A_i))^2$$

$$\leq Cp \sum_i m(B(h(x_i), C_2 diam\ h(A_i))) \leq C\ p\ m(h(\pi_{z_0}^{-1}(B(\sigma, Cr) \cap S(z_0, 1)))).$$

For Lebesgue almost every $\sigma \in \mathbb{S}^1$ there exists $C = C(\sigma)$ so that

$$\frac{m(h(\pi_{z_0}^{-1}(B(\sigma, R) \cap S(z_0, 1))))}{R} \leq C,$$

provided that $R$ is small. Then

$$\sum_i diam\ h(A_i) \leq C\ pr \leq C\ (m_1(F))^{1/2},$$

which finish the proof of the theorem. □

## 6. Rigidity theorem

Let $f$ be a generalized polynomial like map such that $diam\ C_n \to_{n \to \infty} 0$. Let $V_i$, $i \leq N = N(n)$, be the puzzle pieces at level $n$ and let $U_j$, $j \leq N(n+1) - 1$ be the puzzle pieces at level $n+1$ which does not contains the critical point. Define the map $g_n \colon \cup_i U_i \to \cup_j V_j$, as the map $f$ restricted to $\cup_i U_i$. Let $K(g_n)$ be the filled-in Julia set of $g_n$. We can provide each connected component of the Riemann surface $\mathcal{V} = \cup_j V_j$ with its Riemanian hyperbolic metric. The following result is trivial:

**Proposition 6.1.** *There exist $\alpha_2 > \alpha_1 > 1$, with $\alpha_j = \alpha_j(n)$ so that, for any $x \in K(g_n)$, we have*

- $\alpha_1 \leq ||Df(x)||_{\mathcal{V}} \leq \alpha_2$.
- *For any $m > 0$, let $j$ be so that $f^m(x) \in V_j$. Then there is a domain $D(x, m)$ which contains $x$ so that $f^m$ is univalent in $D(x, m)$ and $f^m(D(x, m)) = V_j$.*

The proof of the following corollary is straightforward:

**Corollary 6.1.** *For $n > 0$, there exists a centered Markov family $\mathcal{G}_n$ so that*

- *$\mathcal{G}_n$ covers $K(g_n)$ in all scales in $(C_1, \epsilon_0)$-uniform way, where $(C_1, \epsilon_0)$ depends on $f$ and $n$.*
- *$\mathcal{G}_n$ has $C_2(f, n)$-uniform bounded geometry.*

*Furthermore, this centered Markov family is topologically defined in the following sense: Let $f_j \colon \cup_i U_i^j \to V^j$, $j = 1, 2$, be two generalized polynomial-like maps and let $h \colon \mathbb{C} \to \mathbb{C}$ be a homeomorphism so that*

- *$h(U_i^1) = U_i^2$,*
- *$h(V^1) = V^2$,*
- *$h$ is a conjugacy between $f_1$ and $f_2$.*

*Then $h(\mathcal{G}_n^1) = \mathcal{G}_n^2$.*

*Proof.* Let $x \in K(g_n)$. Then $(A, x) \in \mathcal{G}_n$ if and only if there exists $m \geq 0$ so that $g_n^m(x) \in U_i$ and $A = D(x, m) \cap g_n^{-m}(U_i)$. □

**Corollary 6.2.** *The set $\cup_n K(g_n)$ has zero 2-dimensional Lebesgue measure.*

**Corollary 6.3.** *$K(f)$ has empty interior.*

*Proof.* If $K(f)$ has non empty interior, and by the non existence of wandering components (Sullivan) and the classification of the periodic components, the critical set must be in the interior of $K(f)$, which is a contradiction, since $diam\ C_n \to 0$. □

The following is a Levin and van Strien argument (proposition 3.1 in [LvS]):

**Proposition 6.2.** *Let $f\colon \cup_i U_i \to V$ be a generalized polynomial like map so that*
- *$diam\ C_n \to_n 0$,*
- *There exist $M > 0$ and a subsequence of critical pieces $C_{n_k}$ which admits an annulus $A_k$ so that*
  - *the internal boundary of $A_k$ is $\partial C_{n_k}$,*
  - *$mod\ A_k \geq M$,*
  - *For any $x$ which is not contained in $C_{n_k}$ and $m > 0$ minimal so that $f^m(x) \in C_{n_k}$, there exists a neighborhood $V(x,m)$ of $x$ so that*
    * *$f^m(V(x,m)) = C_{n_k} \cup A_k$,*
    * *$f^m$ is univalent in $V(x,m)$.*

*Then any point $x \in K(f)$ admits a sequence of puzzle pieces which shrinks to $x$. In particular, $K(f)$ is a Cantor set.*

*Proof.* By Proposition 6.1, we can assume that the forward orbit of $x$ accumulate in the critical point. Let $m_k$ be minimal so that $f^{m_k}(x) \in C_{n_k}$. Since $m_k$ is the first entry time, there exists a domain $V(x,k)$ containing $x$ which is mapped univalently by $f^{m_k}$ to $C_{n_k} \cup A_k$. Denote by $P(x,k) \subset V(x,k)$ the piece so that $f^{m_k}P(x,k) = C_{n_k}$. We claim that $diam\ P(x,k)$ goes to zero. Indeed, consider the core curve $\gamma_k$ of the annulus $A_k$. Since $mod\ A_k > M$ and by the Koebe distortion lemma, we have that $\tilde{\gamma}_k = f^{-m_k}(\gamma_k) \cap V(x,k)$ is a curve which contains a disc whose diameter is comparable to $diam\ \tilde{\gamma}_k$ and whose center is $x$. So, if $diam\ \tilde{\gamma}_k$ does not converges to zero, the diameter of this balls also does not converge, and we can find a subsequence so that these balls converge to a ball which is contained in $K(f)$, since $f^m$, $m > 0$, is a normal family in this ball. This is impossible, since $int\ K(f)$ is empty. □

*Proof of Theorem 1.* Since the argument is quite symmetric with respect $f_1$ and $f_2$, we will denote both by $f$. Decompose $K(f)$ in two disjoint sets $K(f) = X_1 \cup X_2$, where $X_2 = \cup_n K(g_n)$. Denote $E_n = K(g_n) - K(g_{n-1})$, with $K(g_0) := \phi$. Then $X_2 = \cup_n E_n$. Define the centered Markov family $\mathcal{F}_2$ in the following way: if $x \in E_n$, then the domains in $\mathcal{F}_2$ with center in $x$ are exactly the domains with center in $x$ which belongs to $\mathcal{G}_n$. It is easy to see that $\mathcal{F}_2$ is topologically defined, covers $X_2$ in all scales and have bounded geometry.

Now we will deal with $X_1$: this is the set of points in $K(f)$ whose positive orbit accumulate in the critical point. Define the Markov family $\mathcal{F}_1$ is the following way: for $x \in K_1$, let $m_k \geq 0$ be minimal so that $f^{m_k}(x) \in C_{n_k}$. By assumption, there is a domain which contains $x$ so that $f^{m_k}$ maps it univalently to $C_{n_k} \cup A_k$. So there exists a puzzle piece $V(x,n)$ which contains $x$ so that
- $f^{m_k}$ is univalent in $V(x,n)$,
- $f^{m_k}(V(x,n)) = C_{n_k}$.

- Moreover, since $mod\ A_k > M$ and by Koebe distortion lemma, there exists $\lambda < 1$ so that for any $(x, n)$ there is $y \in V(x, n)$ satisfying

$$B(y, \lambda\ V(x, n)) \subset V(x, n)$$

where $C = C(M, K_1, K_2)$.

Note that, for a fixed $x$, $diam\ V(x, n) \to_n 0$. Define

$$\mathcal{F}_1 = \{(V(x, n), x)\ :\ x \in X_1\ and\ n > 0\}$$

So $\mathcal{F}_1$ covers $X_1$ in arbitrary small scales, have bounded geometry and, moreover, it is topologically defined. By the Removability Theorem (Theorem 4), $h$ is quasiconformal in $\mathbb{C}$. □

## 7. An illustrative example: Fibonacci combinatorics with degree two

We will prove that

**Theorem 5.** *There exists an unique parameter in the family $x^2 + c$, $c \in \mathbb{R}$, which has Fibonacci combinatorics.*

This theorem is consequence of Yocooz's work (see also Branner [Br], Lyubich & Milnor [LM] and Lyubich [Lyu93]). Since there are easier methods to prove Theorem 5, this is not a "serious" application of Theorem 1. We justify a new approach to this result by its educational interest.

*Proof.* Consider two parameters $c_1$ and $c_2$, so that $P_1(x) = x^2 + c_1$ and $P_2(x) = x^2 + c_2$ have Fibonacci combinatorics. Since these polynomials does not support invariant line fields in its Julia set (for instance, see Levin & van Strien [LvS]), to prove the theorem it is sufficient to prove that there is a quasiconformal conjugacy between these polynomials. We can replace $P_1$ and $P_2$ for induced maps which are generalized polynomial-like maps (whose domains have two connected components) with Fibonacci combinatorics. Since the Julia sets of these generalized polynomial-like maps are totally disconnected and they are topologically conjugated in the real line, by the Sullivan pullback argument we can construct a conjugacy between $P_1$ and $P_2$ which is quasiconformal outside their Julia sets. We can assume that $P_1$ has bounded geometry in the real trace of the principal nest and $P_2$ does not (the other cases are similar). By the starting condition property (see Lyubich & Milnor [LM]), the decay of geometry holds for $P_2$. By Lyubich work [Lyu93] (for $P_2$) and Buff work [Bu] (for $P_1$, see also the previous work by van Strien & Nowicki [vSN94]), there is a sequence of critical puzzle pieces $C_{n_k}$, for both polynomials, which satisfies the hypothesis of Theorem 1 (Indeed, the shape of the critical puzzle pieces in the principal nest of $P_2$ converges to the Julia set of the polynomial $x^2 - 1$, and for $P_1$, for the shape of the Julia set of a polynomial-like map which has an attractor of period two). This finish the proof. □

The interesting point is that we did not use that for $P_1$ and $P_2$ the decay of geometry holds (naturally you can construct a model as in Lyubich & Milnor [LM] which has decay of geometry and, using the above proof, conclude that all real generalized polynomial-like map of degree two with Fibonacci combinatorics have decay of geometry) or that their Julia sets are removable to any quasiconformal map. This method could be very useful to obtain rigidity in higher criticalities, where there are combinatorics where neither the decay of geometry nor the bounded combinatorics holds (See section 17).

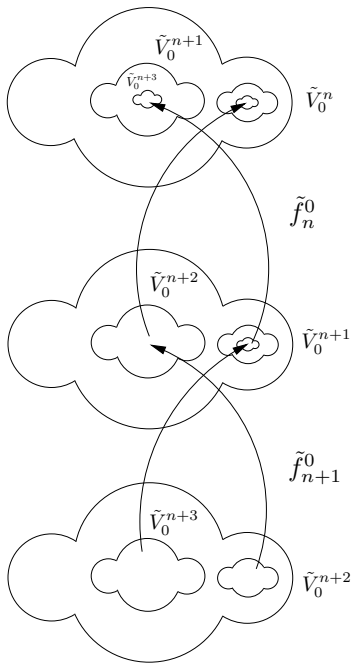

Figure 1.

## Part 2. **Hyperbolicity of the Fibonacci cycle**

### 8. On the maps in the Fibonacci cycle

The most important property of the sequence of critical puzzle pieces in the principal nest (recall the notation in section 4.3) of $\tilde{f}$ is that the shape of them is under control, as we will describe below. An important advise: the figures in this work are only illustrative: the actual shape of the puzzles and Julia sets represented there depends on the criticality $d$ and they are much more complicated. Taking the figures has more actual than they really are could do the reader deduce that $f_n^0$ has degree two, what is not the case. The actual domain $D_0$ represented in Fig. 3, for instance, has more holes (which do not intersect the real line). Realistic pictures of puzzle pieces of Fibonacci maps can be found in Buff[Bu].

See the Fig. 1. It describes how certain puzzles pieces in real Fibonacci maps are permuted by the action of the generalized renormalizations. Buff [Bu] described precisely the shape of $V_0^n$. Denote $d_H(T_1, T_2) := max\{d_h(\overline{T}_1, \overline{T}_2), d_h(\partial T_1, \partial T_2)\}$, where $d_h(\cdot, \cdot)$ represents the Hausdorff metric under the set of compact sets.

**Proposition 8.1** ([Bu]). *Let $\beta$ be the closest period two point to the critical point of $\tilde{f}$. Then there exists a polynomial-like map $g$, with an unique critical point and a hyperbolic attractor of period two, and $0 < \beta < 1$ so that $-1/\beta^2 g^2(\beta z) = g(z)$,*

$$\tilde{f}_0^0(z) = g(\frac{1}{\beta} \cdot z), \ and$$

$$\tilde{f}_0^1(z) = \mp \frac{1}{\beta} g(z).$$

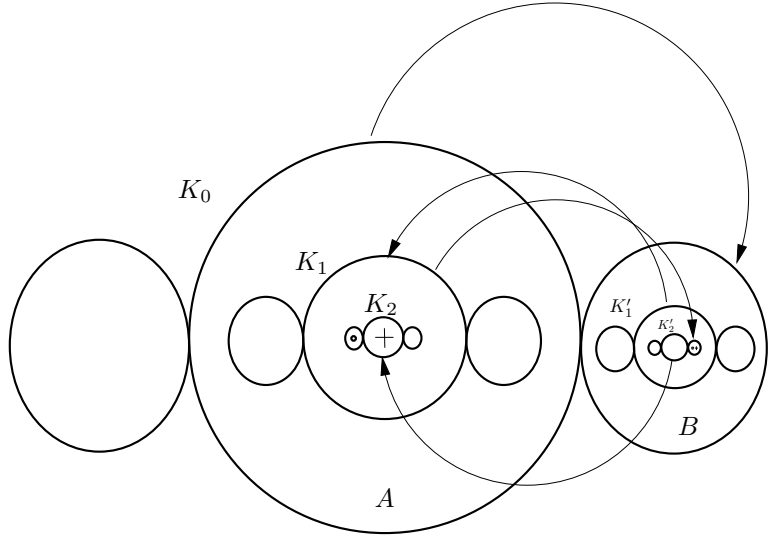

FIGURE 2.

*In particular $g(z) = 1/\beta^n \tilde{f}_n^0(\beta^{n+1} \cdot z)$. Furthermore*

$$dist_H(\frac{1}{\beta^n}\tilde{V}_0^n, K(g)) \to_n 0.$$

Indeed, Buff claimed that $1/\beta^n \overline{\tilde{V}_0^n}$ converges to $K(g)$ in the Hausdorff metric on compact sets, but the slight modification above can be proved in an analogous way (indeed, we will do something similar later, in the proof of proposition 10.1). We will normalize $\tilde{f}$ so that 0 is a maximum and 1 as the fixed point of $\tilde{f}$ in the central domain so that $D\tilde{f}(1) < 0$. Then the correspondent fixed point to $\mathcal{R}\tilde{f}$ is $\beta$ (or $-\beta$).

Denote by $A$ the connected component of the interior of the $K_0 := K(g)$ that contains the critical point and by $B$ the connected component what contains the critical value. We will say that a set $O_1$ is compactly contained in $O_2$, denoted $O_1 \Subset O_2$, if $\overline{O}_1 \subset int\, O_2$.

The following result will be very important:

**Proposition 8.2** ([Bu])**.** *Denote $K_1 := \beta K_0$, $K_2 := \beta^2 K_0$, $K_i' := g^{-1}K_i \cap B$. Then*

- $-1/\beta^2 g^2(\beta z) = g(z)$,
- $K_1 \subset A$,
- $K_2 \Subset K_1$,
- $g(K_1) = K_2'$.

*In particular:*

- $\tilde{f}^0(K_1) = K_0$,
- $\tilde{f}^1(K_1') = K_0$,
- $\tilde{f}^1(K_2') = K_1$,

*Proof.* These items are included or they are an easy consequence of Buff [Bu]. □

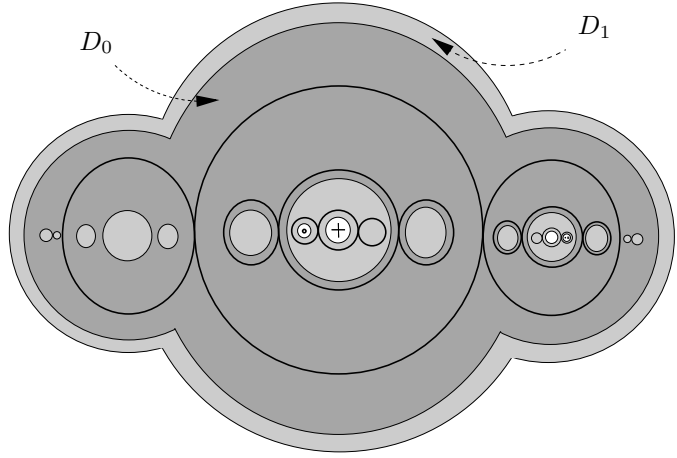

Figure 3.

Fig. 2 describes how $g$ acts on $K_0$, $K_1$, $K_0'$ and $K_1'$: these sets are not represented in the figure, but only the (boundaries of the) connected components (in the figure, circles) of their interior which corresponds to the topological discs $A$, $B$ and $-B$ contained in $K(g)$.

Note that $1/\beta^n \partial \tilde{V}_0^n$ converges to $\partial K_0$, $1/\beta^n \partial \tilde{V}_0^{n+1}$ converges to $\partial K_1$ and $1/\beta^n \partial \tilde{V}_1^{n+1}$ converges to $\partial K_1'$. Since $P(\tilde{f}) \cap (\tilde{V}_0^{n+1} - (\tilde{V}_0^{n+1} \cup \tilde{V}_1^n)) = \phi$, we conclude that $\tilde{f}$ satisfies all the conditions imposed to $f_1$ in theorem 1, with $C_{i_n} := V_0^n$.

The following corollary will be an useful tool in section 11:

**Corollary 8.1.** *There exists an open set $D_1$ which contains $J(g)$ in its interior so that:*

- $g^2 \colon g^{-2}D_1 \to D_1$ *is an unramified covering of degree* $(\deg g)^2$,
- $g^{-2}D_1 \Subset D_1$,
- $\beta g^{-2}D_1 \Subset D_1$,
- $g^{-1}(\beta g^{-2}D_1) \cap B \Subset D_1$.

*Proof.* Let $A_1 \Subset B$ and $B_1 \Subset B$ be domains whose boundaries are Jordan curves very close to $J(g)$ and so that (see fig. 2):

- $B_1 \Subset g^{-1}A_1$,
- $A_1 \Subset g^{-1}B_1$.

Let $U$ be a small neighborhood of $K(f)$ so that $g^{-1}U \Subset U$. We define $D_1$ as (see 8)

$$D_1 := U \setminus (\beta^2 \overline{A}_1 \cup \beta^2 \overline{B}_1 \cup (g^{-1}(\beta^2 \overline{A}_1) \cap B) \cup (g^{-1}(\beta^2 \overline{A}_1) \cap B))$$

Note that

$$D_0 := g^{-2}D_1 \Subset g^{-1}U \setminus (\beta g^{-1}A_1 \cup \beta g^{-1}B_1 \cup g^{-1}(\beta g^{-1}A_1 \cup \beta g^{-1}B_1)) \subset D_1$$

The other properties of $D_1$ are immediate consequence of the previous proposition (provided we have selected $\partial A_1$, $\partial B_1$ and $\partial U$ close of $J(g)$). $\square$

## 9. The Fibonacci renormalization operator: complex analytic version

We are going to define the complex analytic version of the Fibonacci operator, defined in a convenient Banach slice: Let $\tilde{f}$ be one of the maps in the Fibonacci cycle. By proposition 8.2, the $\delta$-neighborhoods of $K_2$, $\tilde{f}(K_2)$ and $\beta K_1'$ are contained in the interior of $K_1 \cup K_1'$, for some $\delta > 0$. Recall that we normalized $\tilde{f}$ so that $\beta$ (or $-\beta$) is a fixed point of $\tilde{f}^2$. Taking $\epsilon_0$ small enough, we can assume that $\beta$ (or $-\beta$) has an analytic continuation, denoted $\beta_f$, for every $f$ so that $|f - \tilde{f}|_{\mathcal{B}(K_1 \cup K_1')} \leq \epsilon_0$. The above considerations implies that there exist $\epsilon_0 > 0$ and $C < 1$ so that if $f \in \mathcal{B}_{nor}(K_1 \cup K_1')$ and $|f - \tilde{f}|_{\mathcal{B}(K_1 \cup K_1')} \leq \epsilon_0$, then $f^2$ is well defined in a $C\delta$-neighborhood of $\beta_f K_1$ and $f$ is defined in a $C\delta$-neighborhood of $\beta_f K_1'$. We will select a domain which contains $K_1 \cup K_1'$ and whose boundary is more regular: Let $U_0$ and $U_1$ be simply connected domains so that $K_1 \Subset U_0 \subset C\delta/2\text{-}K_1$ and $K_1' \Subset U_1 \subset C\delta/2\text{-}K_1'$ and whose boundaries are analytic curves. So for every $f \in \mathcal{B}_{nor}(U_0 \cup U_1)$ satisfying $|f - \tilde{f}|_{\mathcal{B}(U_0 \cup U_1)} \leq \epsilon_0$, $f^2$ is defined in a $2C\delta/3$-neighborhood of $\beta_f \overline{U}_0$ and $f$ is defined in a $2C\delta/3$-neighborhood of $\beta_f \overline{U}_1$. Define the map $\hat{\mathcal{R}} f \in \mathcal{B}_{nor}(2C\delta/3\text{-}\overline{U}_0 \cup 2C\delta/3\text{-}\overline{U}_1)$ as

$$\frac{1}{\beta_f} f^2(\beta_f \cdot z)$$

on $2C\delta/3\text{-}\overline{U}_0$, and

$$\frac{1}{\beta_f} f(\beta_f \cdot z)$$

on $2C\delta/3\text{-}\overline{U}_1$. The **Fibonacci renormalization operator** is defined as $\mathcal{R} := i \circ \hat{\mathcal{R}}$, where $i \colon \mathcal{B}_{nor}(2C\delta/3\text{-}\overline{U}_0 \cup 2C\delta/3\text{-}\overline{U}_1) \to \mathcal{B}_{nor}(U_0 \cup U_1)$ is the inclusion $i(f) = f$. This operator is an example of induced transformation, whose properties will be described in section 13. In the rest of this work, denote $U := U_1 \cup U_2$. Of course, using the complex bounds result to real Fibonacci maps [vSN94], we can define a complex analytic version to a Banach slice $\mathcal{B}_{nor}(\hat{U})$ of an iteration of the Fibonacci renormalization operator in a more easy way, without to use the special features of $\tilde{f}$, and carry out a proof of the hyperbolicity in this Banach space as well. The problem is that exposition would be a little more complicated, specially because we need to deal with an iteration of the Fibonacci operator (not with the operator itself) and there is not warranty in this approach that $K_1 \cup K_1' \subset \hat{U}$: so there is not advantage of to use this approach in the case of the Fibonacci renormalization (but this alternative method is, possibly, the unique way of define the complex extension in the case of generalizations of this work for other combinatorics, as non periodic ones).

## 10. Control of the geometry in the principal nest

Let $f$ be a function in $\mathcal{B}_{nor}(U)$ so that $\mathcal{R}^n f$ exists, for $n \geq 0$, and

$$|\mathcal{R}^n f - \tilde{f}|_{\mathcal{B}(U)} \leq \delta.$$

We will consider a generalized polynomial like restriction to $f$ and prove that the shape of certain puzzle pieces are under control, provided that $\delta$ is small. Let $D \Subset g(D) \subset V$ be simply connected domains, which contains $K(g)$ so that $g \colon D \to g(D)$ is ramified covering map and $D_a \subset A$, $D_b \subset B$, with $\partial D_a \cup \partial D_b$ very close to $J(g)$, simply connected domains so that $g(D_a) \Subset D_b$ and $g(D_b) \Subset D_a$. Define $E_0 := D \setminus \overline{D_a} \cup \overline{D_b}$ and $E_n := g^{-n} E_0$. Then

$$g \colon E_{n+1} \to E_n$$

is a covering map and, because $g$ is a hyperbolic map, if we take a Jordan curve $\gamma \subset E_0$ that is the boundary of bounded domain which contains the postcritical set (and the attractor) of $g$, then $g^{-n}\gamma$ is also Jordan curve which bounded a domain which contains the postcritical set of $g$ and $d_H(g^{-n}\gamma, J(g)) \leq C\lambda^n$, for some $\lambda < 1$ (recall that $d_H(T_1, T_2)$ denotes $max\{d_h(\overline{T}_1, \overline{T}_2), d_h(\partial T_1, \partial T_2)\}$, where $d_h(\cdot, \cdot)$ represents the Hausdorff metric under the set of compact sets). Consider an analytic curve $\gamma \subset E_0$, as above, and denote by $V_0^0$ the simply connected bounded domain whose boundary is $\gamma$. Denote $\tilde{\tau}_n = \beta^n$ and $\beta_i := \beta_{\mathcal{R}^i f}$, and $\tau_n := \beta_{n-1} \cdots \beta_0$. This notation will be used in the rest of this work.

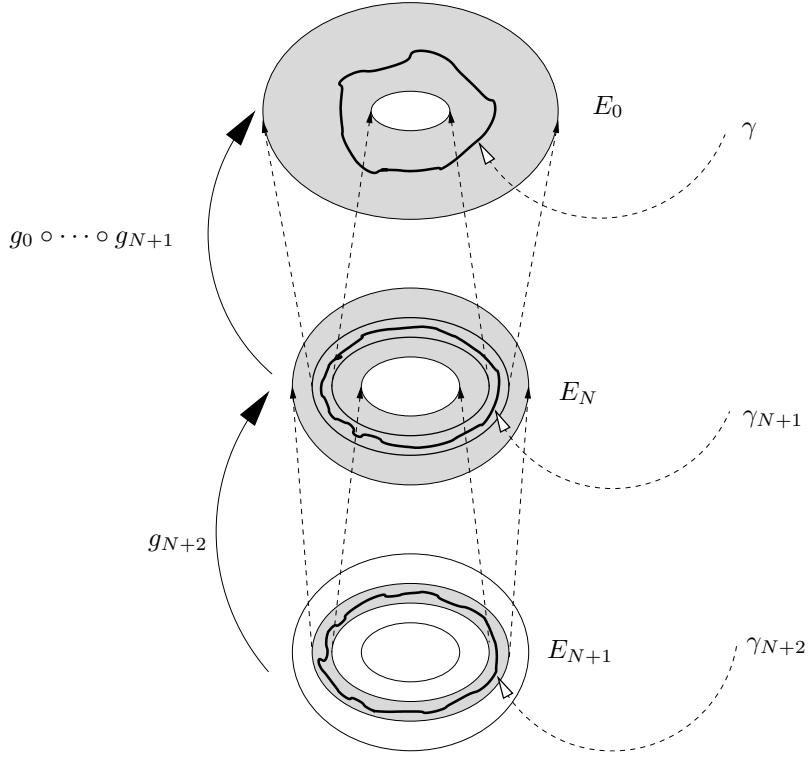

Figure 4.

**Proposition 10.1.** *For every $\epsilon > 0$ there exists $\delta > 0$ and $N = N(\epsilon)$ with the following properties: If $f \in \mathcal{B}_{nor}(U)$ satisfies*

$$|\mathcal{R}^n f - \tilde{f}| \leq \delta,$$

*for $n \geq 0$, and $V_0^1$ and $V_1^1$ are the connected components of $f^{-1} V_0^0$ which contain the critical point and $f(0)$, then $f \colon V_0^1 \cup V_1^1 \to V_0^0$ is a generalized polynomial-like so that*

1. *The critical point does not escape,*

2. *Using the notation defined in section 4.3, $V_0^{n+1} \neq V_1^{n+1}$ and $P(f) \cap V_0^n \subset V_0^{n+1} \cup V_1^{n+1}$.*
3. *$f_n^0 \colon V_0^{n+1} \to V_0^n$ is a restriction of $f^{S_n}$ to $\tau_n U_0$ and $f_n^1 \colon V_1^{n+1} \to V_0^n$ is a restriction of $f^{S_{n-1}}$ to $\tau_n U_1$.*
4. *We have*

$$dist_H(\frac{1}{\tau_n}V_0^n, K_0) \leq \epsilon,$$

$$dist_H(\frac{1}{\tau_n}V_0^{n+1}, K_1) \leq \epsilon,$$

$$dist_H(\frac{1}{\tau_n}V_1^{n+1}, K_1') \leq \epsilon,$$

$$dist_H(\frac{1}{\tau_n}(f_n^1)^{-1}(V_0^{n+1}), K_2') \leq \epsilon,$$

*for $n \geq N$.*

*Furthermore, there exists $C > 0$ so that*

(1) $$dist(\frac{1}{\tau_n}f_n^0(0), \frac{1}{\tau_n}\partial(f_n^1)^{-1}(V_0^{n+1})) \geq C.$$

*In particular, all the assumptions about the generalized polynomial-like maps in Theorem 1 are satisfied by $f$ and, if $\epsilon$ is small enough (see Fig. 5)*

- $V_0^{n-2} \setminus V_0^{n-1} \cup V_1^{n-1} \Subset \tau_{n-2}D_1$,
- $D_0^{n-2} \Subset \tau_{n-2}D_1$,
- $(f_{n-2}^1)^{-1}(\tau_{n-2}D_0^{n-2}) \Subset \tau_{n-2}D_1$,
- $\tau_{n-1}D_0^{n-1} \Subset \tau_{n-2}D_1$.

*where $D_0^n := \tau_{n+1}^{-1}(f_n^0)^{-1}\tau_n D_1$.*

*Proof.* Note that, by definition $\mathcal{R}^n f$ is an appropriated normalization of $(\mathcal{R}^{n-1}f)^2$ in $U_0$ and a normalized restriction of $\mathcal{R}^{n-1}f$ in $U_1$. So it is easy to conclude inductively that

(2) $$f^{S_n}(z) = \tau_n \mathcal{R}^n f(\tau_n^{-1}z) \text{ on } \tau_n U_0,$$

and

(3) $$f^{S_{n-1}}(z) = \tau_n \mathcal{R}^n f(\tau_n^{-1}z) \text{ on } \tau_n U_1.$$

Recall that $\tilde{f}(\beta z) = g(z)$ on $U_0$. Since the renormalizations of $f$ are very close to $\tilde{f}$, $g_n(z) := \mathcal{R}^n f(\beta_n z)$ are very close to $g(z)$ on $U_0$. Consider $N$ so that $dist_H(E_N, J(g)) \leq \epsilon/2$. Choosing $\delta$ small enough, we can assume that:

$$(g_{N+1} \circ \cdots \circ g_0)^{-1}E_0 \Subset E_{N+1},$$

and, for $n \geq 0$,

$$(g_n)^{-1}E_N \Subset E_N.$$

Fig. 4 explains this (but the topology and shape of the sets $E_i$ are not correct in the figure: $E_0$ has two holes, and the number of holes of $E_i$ grows exponentially fast with $i$). Since $\gamma := \partial V_0^0 \subset E_0$, by Eq. 1 and the above properties, there exists Jordan curves $\gamma_n$ defined inductively in the following way: $\gamma_0 := \gamma$, and $\gamma_n$ is the Jordan curve included in the set $f^{-S_{n-1}}\gamma_{n-1}$ so that the bounded domain whose boundary is $\gamma_n$, denoted $\hat{V}_0^n$, contains the critical point. These domains have the following properties:

- $\hat{V}_0^n \Subset \hat{V}_0^{n-1} \subset \tau_{n-1}U_0$,

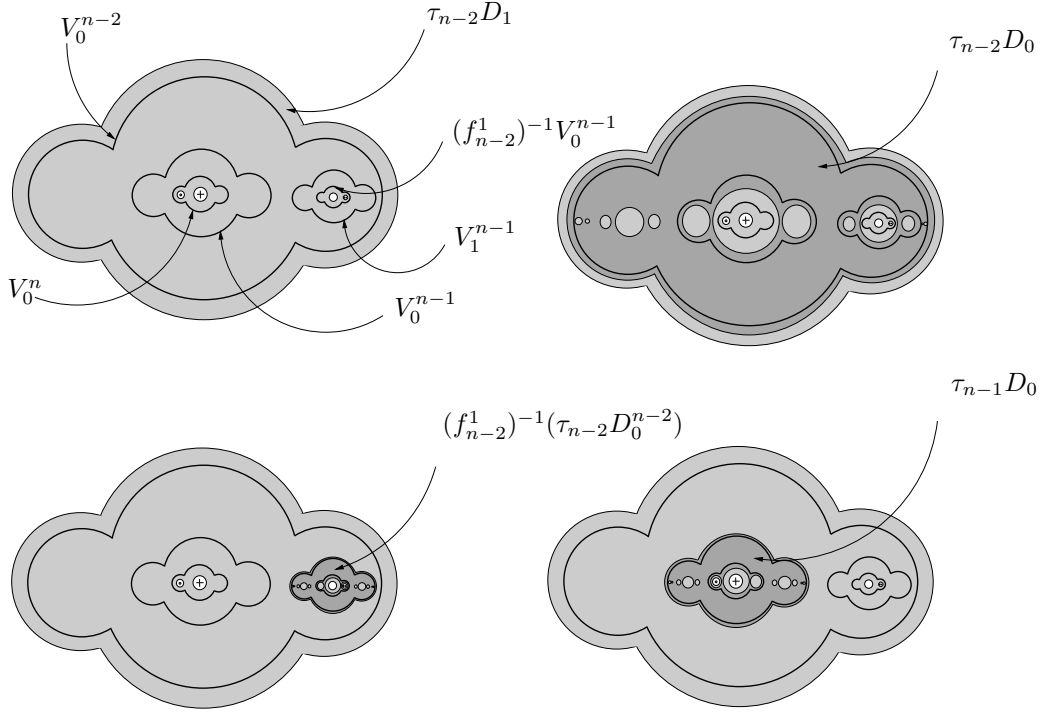

Figure 5.

- $f^{S_n}: \hat{V}_0^{n+1} \to \hat{V}_0^n$ is a ramified covering with an unique critical point.
- We have
$$dist_H(\tau_n^{-1}\hat{V}_0^n, K(g)) \leq \frac{\epsilon}{2},$$
for $n \geq N$. Moreover, the same quantity is smaller than $\epsilon_0$, for all $n \geq 0$.

In particular, the critical point does not escape. We claim that $\hat{V}_0^n = V_0^n$, where $V_0^n$ were defined in section 4.3. Suppose, by induction, that we proved the claim to $V_0^i$, $i < n$. In particular, the time of return of the critical point to $V_0^{n-3}$ and $V_0^{n-2}$ is $S_{n-3}$ and $S_{n-2}$. By definition, $V_0^n$ is the domain of the first return map to $V_0^{n-1}$ which contains the critical point. So it is sufficient to prove that the time of return $m$ of the critical point to $V_0^{n-1}$ is $S_{n-1}$. Clearly $m \leq S_{n-1}$. But note that $f^{S_{n-2}}(0) \in V_0^{n-2} \setminus V_0^{n-1}$, because

(4) $$f^{S_{n-2}}(0) = \tau_{n-2}\mathcal{R}^{n-2}f(0) \in \tau_{n-2}U_1,$$

and
$$V_0^{n-1} \subset \tau_{n-1}U_0 \subset \tau_{n-2}U_0.$$

So the first chance of the forward orbit of $f^{S_{n-2}}(0)$ enter to $V_0^{n-1}$ is its iteration by the time of return of $V_0^{n-2}$ to $V_0^{n-3}$, in other words, $f^{S_{n-2}+S_{n-3}}(0) = f^{S_{n-1}}(0)$. This finish the proof of the claim. Note that, because $f^{S_{n-1}}(0) \in V_0^{n-1} \setminus V_0^n$, we have that $V_1^n \neq V_0^n$, and since $f^{S_{n-1}+S_{n-2}}(0) \in V_0^n$, by Eq. 4 we have $V_1^n = (f_n^1)^{-1}V_0^{n-1}$, where $f_n^1$, as defined in section 4.3, is a restriction of $f^{S_{n-2}}$ on $\tau_{n-1}U_1$. By Eq. 3, and because $\mathcal{R}^n f$ is very close to $\tilde{f}$ and $\tau_n^{-1}V_0^n$ is $\epsilon$-close to $K_0$, we have
$$dist_H(\tau_n^{-1}V_1^{n+1}, K_1') \leq O(\epsilon).$$

Note that the map $f_n$ on $V_0^n$ is the composition of two iterations of $f_{n-1}$: one of them restricted to $V_0^{n-1}$ and another to $V_1^{n-1}$. The map $f_n$ on $V_1^n$ is just one iteration of $f_{n-1}$ on $V_0^{n-1}$. This is sufficient to conclude that the postcritical set of $f$ does not intercept $V_0^i \setminus (V_0^{i+1} \cup V_1^{i+1})$. □

## 11. Quasiconformal rigidity

A map $f \in \mathcal{B}_{nor}(U)$ satisfying the assumptions in proposition 10.1 have a dynamical behavior quite similar to the maps in the Fibonacci cycle: the main result in this section is that there exists a hybrid conjugacy between $f$ and either $\tilde{f}_1$ or $\tilde{f}_2$.

If $A$ is a connected domain which does not contains two points in the complex plane, we can provide it with its hyperbolic metric, denoted $dist(\cdot,\cdot)_D$. Consider $f \in \mathcal{B}_{nor}(U)$ very close to $\tilde{f}$. Let $\tilde{V}_0 = V_0^0$ be the simply connected domain whose boundary is $\gamma$ (see the previous section). Define $V_0^1 := (f_0^0)^{-1} V_0$ and $V_1^1 := (f_0^1)^{-1} V_0$ and similarly $\tilde{V}_0^1$ and $\tilde{V}_1^1$. In this section, we will use the notation introduced in section 4.3. The sets and functions defined there to $f$, can be also defined to $\tilde{f}$, which will be denoted with tilded symbols (for instance, $\tilde{f}_n$, $\tilde{V}_0^n$). Then

**Theorem 6.** *There exists $\delta_0 > 0$ with the following property: Assume that $f \in \mathcal{B}_{nor}(U)$ satisfies*
$$|\mathcal{R}^n f - \tilde{f}|_{\mathcal{B}(U)} \leq \delta \leq \delta_0$$
*for $n \geq 0$. Then there exists a $O(\delta)$-quasiconformal map $\phi \colon \mathbb{C} \to \mathbb{C}$ which is a hybrid conjugacy between $f \colon V_0^1 \cup V_1^1 \to V_0$ and $\tilde{f} \colon \tilde{V}_0^1 \cup \tilde{V}_1^1 \to \tilde{V}_0$, with $\phi(V_0) = \tilde{V}_0$.*

As the first step to proof the main result of this section, we will prove:

**Proposition 11.1.** *There exists $\delta_0 > 0$ with the following property: Assume that $f \in \mathcal{B}_{nor}(U)$ satisfies*
$$|\mathcal{R}^n f - \tilde{f}|_{\mathcal{B}(U)} \leq \delta \leq \delta_0$$
*for $n \geq 0$. Then there exists a topological conjugacy $\phi \colon \mathbb{C} \to \mathbb{C}$ between $f \colon V_0^1 \cup V_1^1 \to V_0$ and $\tilde{f} \colon \tilde{V}_0^1 \cup \tilde{V}_1^1 \to \tilde{V}_0$, so that*
- $\phi(V_0) = \tilde{V}_0$,
- $\phi$ is $O(\delta)$-quasiconformal outside $K(f)$.

**Lemma 11.1.** *There exists $\delta_0 > 0$, $\epsilon_0 > 0$ and $C > 0$ with the following property: if*
$$|\mathcal{R}^n f - \tilde{f}|_{\mathcal{B}(U)} \leq \delta \leq \delta_0$$
*for $n \geq 0$ and $\phi_0 \colon A_0 \to \tilde{A}_0$, where $A_0 = \overline{V_0^0} \setminus \cup V_i^1$, is an injective continuous maps so that*
- $|\phi_0(z) - z| \leq \epsilon$, for $x \in A_0$, for some $\epsilon$ so that $\epsilon_0 > \epsilon > C\delta$.
- $\phi_0 \circ f = \tilde{f} \circ \phi_0$ on $\partial V_0^1 \cup \partial V_1^1$.

*then there exists a sequence of domains $A_n$, the correspondent tilded domains $\tilde{A}_n$, and maps $\phi_n \colon A_n \to \tilde{A}_n$, $\psi_n \colon A_{n-1} \cup (f_{n-1}^1)^{-1} A_{n-1} \to \tilde{A}_{n-1} \cup (\tilde{f}_{n-1}^1)^{-1} \tilde{A}_{n-1}$ so that*

1. $A_n := (f_{n-1}^0)^{-1}(A_{n-1} \cup (f_{n-1}^1)^{-1}(A_{n-1}))$.
2. $A_n \subset \overline{V_0^n} \setminus (V_0^{n+1} \cup V_1^{n+1})$ and $\partial V_0^n$, $\partial V_0^{n+1}$, $\partial V_1^{n+1}$ are connected components of $\partial A_n$.

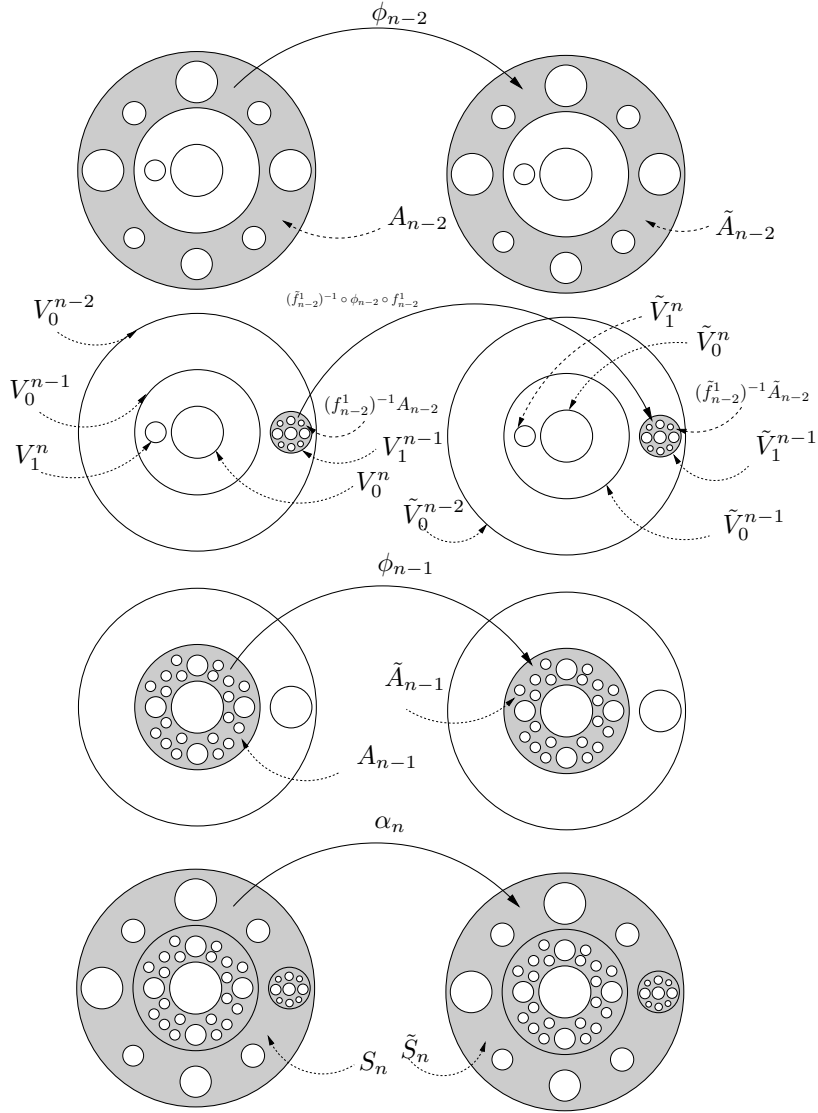

Figure 6.

3. $\psi_n$ is equal to $\phi_{n-1}$ on $A_{n-1}$ and it coincides with $\tilde{f}_{n-1}^{-1} \circ \phi_{n-1} \circ f_{n-1}^1$ on $(f_{n-1}^1)^{-1} A_{n-1}$.
4. $\psi_n \circ f_n^0 = \tilde{f}_n^0 \circ \phi_n$ on $A_n$ and the maps $\phi_{n-1}$, $\psi_n$ and $\phi_n$ coincides on $\partial V_0^n$.
5. for $z \in \tau_n^{-1} A_n$,
$$dist(z, \tilde{\tau}_n^{-1}\phi_n(\tau_n \cdot z))_{D_0} \leq \epsilon.$$
6. for $z \in \tau_{n-1}^{-1} A_{n-1} \cup (f_{n-1}^1)^{-1} A_{n-1}$ we have
$$dist(z, \tilde{\tau}_{n-1}^{-1}\psi_{n-1}(\tau_{n-1} \cdot z))_{D_0} \leq C\epsilon.$$

*Proof.* Recall that we assume that all renormalizations of $f$ and $\tilde{f}$ are very close. Define $\psi_1$ equals to $\phi_0$ on $A_0$ and equals to $(\tilde{f}_0^1)^{-1} \circ \phi_0 \circ f_0^1$ on $(f_0^1)^{-1}(A_0)$. Define $\phi_1$ as the unique map so that

- $\phi_1 = \phi_0$ in $\partial V_1^0$,
- $\psi_1 \circ f_1^0 = \tilde{f}_1^0 \circ \phi_1$ on $A_1$.

$\phi_1$ is well defined if $\epsilon_0 + \delta_0$ is small enough. Since the first renormalizations of $f$ and $\tilde{f}$ are very close, we can conclude that $\phi_1$ and the extension of $\phi_0$ are $O(\delta_0 + \epsilon_0)$-close of identity. Assume by induction that we defined $\phi_i$, $i < n$, with the required properties. Here $D_0$ is as in corollary 8.1. Define $\psi_{n-1}$ equals to $\phi_{n-2}$ on $A_{n-2}$ and as $(\tilde{f}_{n-2}^1)^{-1} \circ \phi_{n-2} \circ f_{n-2}^1$ on $(f_{n-2}^1)^{-1} A_{n-2}$. Denote $S_n := A_{n-2} \cup A_{n-1} \cup (f_{n-2}^1)^{-1} A_{n-2}$. Let $\alpha_n \colon S_n \to \tilde{S}_n$ be the map which is equal to $\psi_{n-1}$ on $A_{n-2} \cup (f_{n-2}^1)^{-1} A_{n-2}$ and it is equal to $\phi_{n-1}$ on $A_{n-1}$ (Fig. 6 describes how $\alpha_n$ is defined in each part of $S_n$). Then define $\phi_n \colon A_n \to \tilde{A}_n$ as the unique map (see Fig. 7) so that

- $\alpha_n \circ f_{n-2}^0 \circ f_{n-1}^0 = \tilde{f}_{n-2}^0 \circ \tilde{f}_{n-1}^0 \circ \phi_n$ on $A_n$, and
- $\phi_n = \phi_{n-1}$ on $\partial V_0^n$.

Note that $\phi_n$ exists because $\tau_i^{-1} f_i^0(\tau_i z)$ is $\delta$-close to $\tilde{f}^0(z) = \tilde{\tau}_i^{-1} \tilde{f}_i^0(\tilde{\tau}_i z)$, and $\tilde{\tau}_{n-2}^{-1} \alpha_n(\tau_{n-2} z)$ is $O(\epsilon + \delta)$-close to identity. Recall that

$$g(z) := \tilde{\tau}_i^{-1} \tilde{f}_i^0(\tilde{\tau}_{i+1} z)$$

is a polynomial-like map with a period two attractor, which does not depends on $i$. In particular, $\tilde{f}_{n-2}^0 \circ \tilde{f}_{n-1}^0(z) = \tilde{\tau}_{n-2} g^2(\tilde{\tau}_n^{-1} z)$. Denote

$$g_i(z) := \tau_i^{-1} f_i^0(\tau_{i+1} z).$$

Since all renormalizations of $f$ are $\delta$-close to $\tilde{f}$, corollary 8.1 remains valid replacing $g^2$ by $g_{n-2} \circ g_{n-1} = \tau_{n-2}^{-1} f_{n-2}^0 \circ f_{n-1}^0(\tau_n z)$. We claim that

$$dist(z, \tilde{\tau}_{n-2}^{-1} \alpha_n(\tau_{n-2} z))_{D_1} \leq \frac{\epsilon}{\lambda} + C\delta, \text{ on } S_n,$$

where $\lambda > 1$ does not depend on $n$. We will divide $S_n$ in three parts, $A_{n-2}$, $A_{n-1}$ and $(f_{n-2}^1)^{-1} A_{n-2}$, and prove of the claim in each one of them. The first part is $A_{n-2}$. Since $D_0 \Subset D_1$ there exists $\lambda_1 > 1$, so that

$$dist(x, y)_{D_1} \leq \frac{1}{\lambda_1} dist(x, y)_{D_0}, \text{ if } x, y \in D_0.$$

Because $\alpha_n$ is equal to $\phi_{n-2}$ on $A_{n-2}$ and the inductive assumption, we get

$$dist(z, \tilde{\tau}_{n-2}^{-1} \alpha_n(\tau_{n-2} z))_{D_1} \leq \epsilon/\lambda_1$$

on $A_{n-2}$.

The second part is $A_{n-1}$. Recall that $\beta D_0 \Subset D_1$, which implies that there exists $\lambda_2 > 1$ so that

$$dist(x, y)_{D_1} \leq \frac{1}{\lambda_2} dist(x, y)_{\beta D_0}, \text{ if } x, y \in \beta D_0.$$

By induction assumption, we have

$$dist(\frac{z}{\beta_{n-1}}, \tilde{\tau}_{n-1}^{-1} \phi_{n-1}(\tau_{n-1} \frac{z}{\beta_{n-1}}))_{D_0} \leq \epsilon$$

on $\tau_{n-2} A_n$. Moreover, since $\beta_n - \beta = O(\delta)$, $A_{n-1} \subset \overline{V}_0^{n-1} \setminus (V_0^n \cup V_1^n)$ and

$$\tau_{n-1}^{-1}(\overline{V}_0^{n-1} \setminus (V_0^n \cup V_1^n)) \sim K_0 \setminus (int\ K_1 \cup int\ K_1') \Subset D_0,$$

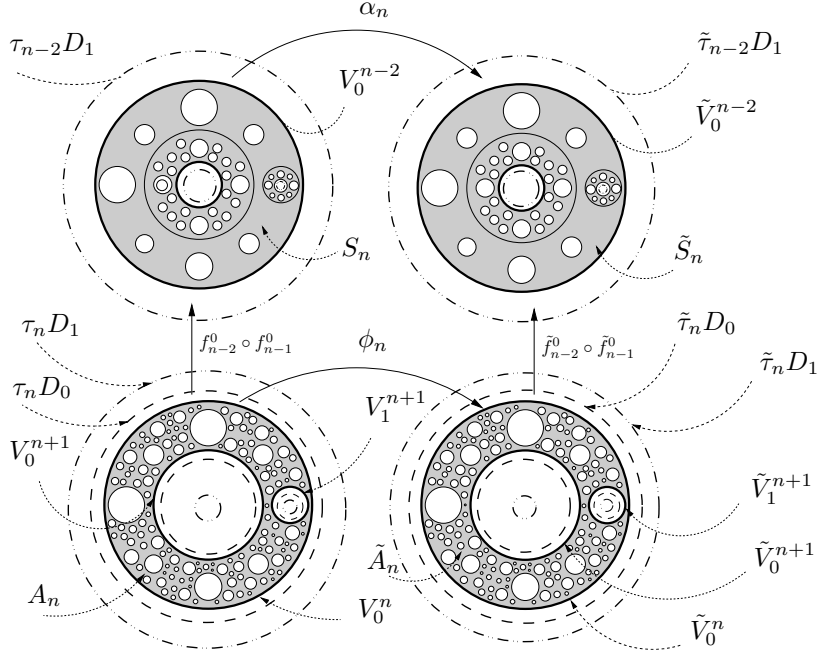

Figure 7.

there exists $C > 0$, which does not depends on $n$, so that, for $z \in \tau_{n-2}A_{n-1}$, we have
$$dist(\frac{z}{\beta}, \frac{z}{\beta_{n-1}})_{D_0} \leq C\delta.$$

So we have
$$dist(z, \tilde{\tau}_{n-2}^{-1}\alpha_n(\tau_{n-2}z))_{D_1}$$
$$= dist(\beta\frac{z}{\beta}, \beta\tilde{\tau}_{n-1}^{-1}\phi_{n-1}(\tau_{n-1}\frac{z}{\beta_{n-1}}))_{D_1}$$
$$\leq \frac{1}{\lambda_2}dist(\beta\frac{z}{\beta}, \beta\tilde{\tau}_{n-1}^{-1}\phi_{n-1}(\tau_{n-1}\frac{z}{\beta_{n-1}}))_{\beta D_0}$$
$$= \frac{1}{\lambda_2}dist(\frac{z}{\beta}, \tilde{\tau}_{n-1}^{-1}\phi_{n-1}(\tau_{n-1}\frac{z}{\beta_{n-1}}))_{D_0}$$
$$\leq \frac{1}{\lambda_2}dist(\frac{\beta_{n-1}}{\beta}\frac{z}{\beta_{n-1}}, \frac{z}{\beta_{n-1}})_{D_0} + \frac{1}{\lambda_2}dist(\frac{z}{\beta_{n-1}}, \tilde{\tau}_{n-1}^{-1}\phi_{n-1}(\tau_{n-1}\frac{z}{\beta_{n-1}}))_{D_0}$$
$$\leq \frac{\epsilon}{\lambda_2} + C\delta$$

on $\tau_{n-2}A_{n-1}$.

The last part is $(f_{n-2}^1)^{-1}A_{n-2}$. Because $\tau_{n-2}^{-1}f_{n-2}^1(\tau_{n-2}z)$ is very close to $\tilde{f}^1 = \tilde{\tau}_{n-2}^{-1}\tilde{f}_{n-2}^1(\tilde{\tau}_{n-2}z)$, we have $\tau_{n-2}^{-1}(f_{n-2}^1)^{-1}\tau_{n-2}D_1 \sim (\tilde{f}_0^1)^{-1}D_1 \Subset D_0$, so there exists $\lambda_3 > 1$ so that
$$dist(x,y)_{\tau_{n-2}^{-1}(f_{n-2}^1)^{-1}\tau_{n-2}D_1} \leq \frac{1}{\lambda_3}dist(x,y)_{D_0},$$

for $x, y \in \tau_{n-2}^{-1}(f_{n-2}^1)^{-1}\tau_{n-2}D_1$. Here $\lambda_3$ does not depends on $n$. Since $\tilde{f}_{n-2}^1 \circ \alpha_n(z) := \phi_{n-2} \circ f_{n-2}^1(z)$ on $(f_{n-2}^1)^{-1}A_{n-2}$, we get

$$dist(\tau_{n-2}^{-1}f_{n-2}^1(\tau_{n-2}z), \tau_{n-2}^{-1}f_{n-2}^1(\tau_{n-2}\tilde{\tau}_{n-2}^{-1}\alpha_n(\tau_{n-2}z)))_{D_1}$$
$$\leq dist(\tau_{n-2}^{-1}f_{n-2}^1(\tau_{n-2}z), \tilde{\tau}_{n-2}^{-1}\tilde{f}_{n-2}^1(\tilde{\tau}_{n-2}\tilde{\tau}_{n-2}^{-1}\alpha_n(\tau_{n-2}z)))_{D_1}$$
$$+dist(\tilde{\tau}_{n-2}^{-1}\tilde{f}_{n-2}^1(\tilde{\tau}_{n-2}\tilde{\tau}_{n-2}^{-1}\alpha_n(\tau_{n-2}z)), \tau_{n-2}^{-1}f_{n-2}^1(\tau_{n-2}\tilde{\tau}_{n-2}^{-1}\alpha_n(\tau_{n-2}z)))_{D_1}$$
$$\leq dist(\tau_{n-2}^{-1}f_{n-2}^1(\tau_{n-2}z), \tilde{\tau}_{n-2}^{-1}\phi_{n-2}(\tau_{n-2}\tau_{n-2}^{-1}f_{n-2}^1(\tau_{n-2}z)))_{D_1}$$
$$+dist(\tilde{\tau}_{n-2}^{-1}\tilde{f}_{n-2}^1(\tilde{\tau}_{n-2}\tilde{\tau}_{n-2}^{-1}\alpha_n(\tau_{n-2}z)), \tau_{n-2}^{-1}f_{n-2}^1(\tau_{n-2}\tilde{\tau}_{n-2}^{-1}\alpha_n(\tau_{n-2}z)))_{D_1}$$
$$\leq \epsilon + C\delta$$

on $\tau_{n-2}^{-1}(f_{n-2}^1)^{-1}A_{n-2}$. Because the map $\tau_{n-2}^{-1}f_{n-2}^1(\tau_{n-2}\cdot)\colon \tau_{n-2}^{-1}(f_{n-2}^1)^{-1}\tau_{n-2}D_1 \to D_1$ is univalent, it preserves the hyperbolic metric, so we get

$$dist(z, \tilde{\tau}_{n-2}^{-1}\alpha_n(\tau_{n-2}z))_{D_1} \leq \frac{\epsilon}{\lambda_3} + C\delta,$$

on $\tau_{n-2}^{-1}(f_{n-2}^1)^{-1}A_{n-2}$. This finish the proof of the claim (take $\lambda := \min\{\lambda_1, \lambda_2, \lambda_3\}$). Note that,

$$dist(\tau_{n-2}^{-1}f_{n-2}^0 \circ f_{n-1}^0(\tau_n z), \tilde{\tau}_{n-2}^{-1}\tilde{f}_{n-2}^0 \circ \tilde{f}_{n-1}^0(\tau_n^{-1}z))_{D_1} \leq C\delta,$$

on an open set what contains $\tau_n^{-1}A_n$ and $\tilde{\tau}_n^{-1}\tilde{A}_n$. So

$$dist(\tilde{\tau}_{n-2}^{-1}\tilde{f}_{n-2}^0 \circ \tilde{f}_{n-1}^0(\tilde{\tau}_n z), \tilde{\tau}_{n-2}^{-1}\alpha_n(f_{n-2}^0 \circ f_{n-1}^0(\tau_n z)))_{D_0}$$
$$\leq dist(\tilde{\tau}_{n-2}^{-1}\tilde{f}_{n-2}^0 \circ \tilde{f}_{n-1}^0(\tilde{\tau}_n z), \tau_{n-2}^{-1}f_{n-2}^0 \circ f_{n-1}^0(\tau_n z))_{D_0}$$
$$+dist(\tau_{n-2}^{-1}f_{n-2}^0 \circ f_{n-1}^0(\tau_n z), \tilde{\tau}_{n-2}^{-1}\alpha_n(\tau_{n-2}\tau_{n-2}^{-1}f_{n-2}^0 \circ f_{n-1}^0(\tau_n z)))_{D_0}$$
$$\leq \frac{\epsilon}{\lambda} + C\delta$$

on $\tau_n^{-1}A_n$. This implies, because $g^2 = \tilde{\tau}_{n-2}^{-1}\tilde{f}_{n-2}^0 \circ \tilde{f}_{n-1}^0(\tilde{\tau}_n\cdot)\colon D_0 \to D_1$ is covering map (hence it preserves the Riemanian hyperbolic metric), that, for each $z \in \tau_n^{-1}A_n$, there exists

$$z_0 \in (\tilde{\tau}_{n-2}^{-1}\tilde{f}_{n-2}^0 \circ \tilde{f}_{n-1}^0(\tilde{\tau}_n\cdot))^{-1}(\tilde{\tau}_{n-2}^{-1}\alpha_n \circ f_{n-2}^0 \circ f_{n-1}^0(\tau_n z))$$

so that $dist(z_0, z)_{D_0} < \epsilon/\lambda + C\delta$. Moreover, since $g^2\colon D_0 \to D_1$ is a covering map, there exists $C_1 > 0$ so that

$$min\{dist(z_1, z_2)_{D_0}\colon z_1, z_2 \in g^{-2}(z), z_1 \neq z_2, z \in \tilde{\tau}_{n-2}^{-1}\tilde{S}_n\} > C_1,$$

where $C_1$ does not depend on $n$. If we assume that $\epsilon/\lambda + C\delta < \epsilon << \alpha/3$, we conclude that, for $z_1 \in (\tilde{\tau}_{n-2}^{-1}\tilde{f}_{n-2}^0 \circ \tilde{f}_{n-1}^0(\tilde{\tau}_n\ \cdot))^{-1}(\tilde{\tau}_{n-2}^{-1}\alpha_n \circ f_{n-2}^0 \circ f_{n-1}^0(\tau_n z))$, $z_1 \neq z_0$, we have $dist(z_1, z)_{D_0} \geq 2C_1/3$. Since

$$\tilde{\tau}_n^{-1}\phi_n(\tau_n z) \in (\tilde{\tau}_{n-2}^{-1}\tilde{f}_{n-2}^0 \circ \tilde{f}_{n-1}^0(\tilde{\tau}_n\cdot))^{-1}(\tilde{\tau}_{n-2}^{-1}\alpha_n \circ f_{n-2}^0 \circ f_{n-1}^0(\tau_n z)),$$

and $A_n$ is connected then either

$$dist(z, \tilde{\tau}_n^{-1}\phi_n(\tau_n z))_{D_0} \leq \epsilon \text{ on } A_n, \text{ or}$$
$$dist(z, \tilde{\tau}_n^{-1}\phi_n(\tau_n z))_{D_0} \geq 2C_1/3 \text{ on } A_n.$$

Since $\phi_n$ coincides with $\phi_{n-1}$ on $\partial V_0^n$ we have $dist(z, \tilde{\tau}_n\phi_n(\tau_n z))_{D_0} < O(\epsilon) << 2C_1/3$, for $z \in \tau_n^{-1}\partial V_n^0$, due to inductive assumption. We finished the proof. □

If $\mathcal{F}$ is a family of domains (as $\mathcal{D}_n$ and $\mathcal{E}_n$), denote $\cup \mathcal{F} := \cup_{D \in \mathcal{F}} D$, $\mathcal{F} \cap D_0 := \{D \in \mathcal{F}\colon D \subset D_0\}$. The following lemma is easy:

**Lemma 11.2.** *Let $f\colon \cup_i V_i^1 \to V_0^0$ be a generalized polynomial-like maps whose critical point does not escape and have a sequence of critical puzzle pieces whose diameter converges to zero. Then there exists a sequence $c_n > 0$ so that, if $D_0 \in \mathcal{D}_n \cap V_0^n$ then there exists an annulus $A$ so that*

- *The internal boundary of $A$ is $\partial D_0$.*
- *$A \subset V_0^n \setminus \overline{\cup \mathcal{D}_n}$.*
- *$\mathrm{mod}\ A \geq c_n$.*

*Proof of Proposition 11.1.* We will use tilded symbols to denote the corresponding objects relatives to $\tilde{f}$. Using holomorphic motions, it is easy to construct an $K$-quasiconformal map $\phi_0$ in the complex plane, with $K = O(\delta)$, so that $\phi_0(V_0^0) = \tilde{V}_0^0$ and $\phi_0 \circ f = \tilde{f} \circ \phi_0$ on $\partial V_0^1 \cup \partial V_1^1$. Apply lemma 11.1 to $\phi_0$. We will prove by induction that there exist $K > 0$ and a sequence of injective continuous maps $\phi_n\colon V_0^0 \to V_0^0$ so that

- $\phi_n$ coincides on $A_n$ with the map $\phi_n$ given by lemma 11.1.
- $\phi_n$ maps the set of points outside the domain of the first return map of $V_0^n$ to the correspondent set to $\tilde{f}$. Furthermore

$$\phi_n \circ f = \tilde{f} \circ \phi_n$$

  in this set.
- $\phi_n$ is $K$-quasiconformal in a neighborhood of every point $x$ which satisfies all the following properties:
  - $x$ does not belong to the domain of the first return map of $V_0^n$.
  - $x$ is not in the Julia set.

  Note that the points in the boundaries of the domains in $\mathcal{D}_n$ satisfies these properties.
- $\phi_n$ is equal to $\phi_{n-1}$ outside the domains of the first return map of $V_0^{n-1}$.

The pointwise limit of $\phi_n$ has a continuous extension to the Julia set of $f$. Moreover, it is a conjugacy between $f$ and $\tilde{f}$ which is $K$-quasiconformal outside $K(f)$.

Assume by induction that we have constructed $\phi_{n-1}$. Let $\{D_j\}$, with $0 \in D_0$, be an enumeration of the subset of domains in $\mathcal{D}_{n-1}$ which are contained in $V_0^{n-1}$. Enumerate the corresponding family for $\tilde{f}$ in the following way: $\tilde{D}_j := \phi_{n-1} D_j$. . Define inductively the family of domains

$$\{D_{i_1,\ldots,i_k}\colon\ (i_1, : , i_k) \in \mathbb{N}_*^k, \text{ for some } k\}$$

($\mathbb{N}_* := \mathbb{N} - \{0\}$) in the following way:

$$D_{i_1,\ldots,i_k} := f_{n-1}^{-k}(D_{i_k}) \cap D_{i_1,\ldots,i_{k-1}}.$$

Note that

(5) $$\mathrm{mod}\ D_{i_1,\ldots,i_{k-1}} \setminus D_{i_1,\ldots,i_{k-1},i_k} \geq c_{n-1},$$

where $c_{n-1}$ is given in lemma 11.2. Define $\pi_{i_1,\ldots,i_k}\colon D_{i_1,\ldots,i_k} \to V_0^{n-1}$ by $f_{n-1}^k$. We are going to define inductively a sequence of continuous injective maps $\psi_n^k\colon V_0^{n-1} \to V_0^{n-1}$ in the following way: $\psi_n^0 := \phi_{n-1}$. We define $\psi_n^k$ as equal to $\psi_n^{k-1}$ outside $\cup_{(i_1,\ldots,i_k)\in\mathbb{N}_*^k} D_{i_1,\ldots,i_k}$ and as equal to

$$\tilde{\pi}_{i_1,\ldots,i_k}^{-1} \circ \psi_n \circ \pi_{i_1,\ldots,i_k}$$

inside $D_{i_1,\ldots,i_k}$. Note that $\psi_n^k$ are continuous and injective. We claim that $\psi_n^k$ converges to some injective continuous map $\psi_n^\infty \colon V_0^{n-1} \to \tilde{V}_0^{n-1}$. To prove this, we will decompose $V_0^{n-1}$ in three parts, $V_0^{n-1} = P_1 \cup P_2 \cup P_3$.

The set $P_1$ consists of the points in $V_0^{n-1}$ whose forward orbit enter in $V_0^n$ at least once. In other words, $P_1 = \cup(\mathcal{E}_n \cap V_0^{n-1})$. The sequence $\psi_n^k$ is eventually locally constant in $P_1$ (eventually locally constant at a point $z$ means that there exists an open set $B$, which contains $z$, so that $\psi_n^i = \psi_n^{i+i}$ on $B$, for large $i$).

The set $P_2$ consists of the points in $V_0^{n-1}$ whose forward orbit returns to $V_0^{n-1}$ only a finite number of times and never enters in $V_0^n$. Note that this set includes the boundaries of domains in $\mathcal{E}_n$ which are contained in $V_0^{n-1}$. In $P_2$ the sequence $\psi_n^k$ is also eventually locally constant. Moreover, if a point $z$ belongs to $P_2 \setminus K(f)$, then we can conclude that the maps $\psi_n^k$ are $K$-quasiconformal in a neighborhood of this point, since the map $\psi_n^k$ is defined as a lift of the map $\phi_n$ by univalent maps, and, by induction hypothesis, $\phi_n$ is $K$-quasiconformal in a neighborhood of an large iteration of $z$.

The last set is $P_3$. This is the set of points in $V_0^{n-1}$ whose forward orbit returns to $V_0^{n-1}$ at infinite many times, but never enters in $V_0^n$. Note that $P_3 \subset K(f)$. We will give a description of $P_3$: define $p_{i_1,i_2,\ldots}$, with $(i_1,\ldots) \in \mathbb{N}_*^{\mathbb{N}}$, as the unique point in
$$\cap_k D_{i_1,\ldots,i_k}.$$
This set contains an unique point due to Eq. 5. Then $P_3 = \{p_\sigma \colon \sigma \in \mathbb{N}_*^{\mathbb{N}}\}$. Note that

(6) $$\psi_n^\ell(D_{i_1,\ldots,i_k}) \subset \tilde{D}_{i_1,\ldots,i_k}, \text{ for } \ell \text{ large enough.}$$

So we have, again due to Eq. 5, that $\psi_n^k(p_{i_1,i_2,\ldots}) \to_k \tilde{p}_{i_1,i_2,\ldots}$.

So the pointwise limit of $\psi_n^k$, denoted by $\psi_n^\infty$, is well defined to every point in $V_0^{n-1}$. The continuity at the points in $P_1 \cup P_2$ is obvious. The continuity in the points in $P_3$ is consequence of Eq. 5 and Eq. 6. Extend the map $\psi_n^\infty$ outside $V_0^{n-1}$ in the following way: if $D \in \mathcal{E}$, $D \neq V_0^{n-1}$, define
$$\psi_{n/D}^\infty := \tilde{\pi}_{n/D}^{-1} \circ \psi_n^\infty \circ \pi_n.$$
Outside $\mathcal{E}_{n-1}$ define $\psi_n^\infty := \phi_{n-1}$. Note that $\psi_n^\infty$ is continuous and injective and satisfies the following properties:

- $\psi_n^\infty$ maps the set of points outside the domain of the first entry map to $V_0^n$ to the correspondent set to $\tilde{f}$ and
$$\psi_n^\infty \circ f = \tilde{f} \circ \psi_n^\infty$$
in this set.
- It is $K$-quasiconformal in an neighborhood of every point $z$ which satisfies the both properties below:
  - $z$ does not belong to the domain of the first entry map to $V_0^n$,
  - $z$ does not belong to $K(f)$.

Suppose that $f_{n-1}^0(0) \in D \in \mathcal{E}_n$. Since the maps $f_{n-1}^0$ and $f_{n-1}^1$ are very close to $\tilde{f}_{n-1}^0$ and $\tilde{f}_{n-1}^1$, and $V_0^{n-1}$, $V_1^{n-1}$ and their boundaries are very close of $\tilde{V}_0^{n-1}$, $\tilde{V}_1^{n-1}$ and their boundaries, the first entry of $f_{n-1}^0(0)$ to $V_0^n$ occurs after apply $f_{n-1}^1$ once. In particular, $D$ is $(f_{n-1}^0)^{-1}V_0^n$, which is very close to $(\tilde{f}_{n-1}^0)^{-1}\tilde{V}_0^n$. Note that $\psi_n^\infty$ is equal to $\psi_n$ on $\partial(f_{n-1}^0)^{-1}V_0^n$, where $\psi_n$ is the map defined in lemma 11.1. In particular

$$dist(z, \tilde{\tau}_{n-1}^{-1}\psi_n^\infty(\tau_{n-1}z)) \leq O(\delta),$$

for $z \in \partial(f_{n-1}^0)^{-1}V_0^n$, which implies that $\tilde{f}_n^0(0) \in \psi_n^\infty((f_{n-1}^0)^{-1}V_0^n)$ (use Eq. 1 in proposition 10.1). Since $\tilde{f}_n^0(0) \in (\tilde{f}_{n-1}^0)^{-1}\tilde{V}_0^n$, we have $\psi_n^\infty((f_{n-1}^0)^{-1}V_0^n) = (\tilde{f}_{n-1}^0)^{-1}\tilde{V}_0^n$. Modify $\psi_n^\infty$ only inside $V_1^{n+1}$ to get a new map $\psi_n$ so that $\psi_n^\infty(f_{n-1}^0(0)) = \tilde{f}_{n-1}^0(0)$, but keeping the $K$-quasiconformality of the map $\psi_n$ in a neighborhood of $\partial D$. Denote by $\{D_j\} \subset \mathcal{D}_n$, with $0 \in D_0$, the family of connected components of the preimages of the elements of $\mathcal{E}_n \cap V_0^{n-1}$ by the map $f_n^0 \colon V_0^n \to V_0^{n-1}$. Define $\phi_n \colon V_0^n \to \tilde{V}_0^n$ as the unique map so that:

- $\phi_n$ is equal to $\psi_n$ in $\partial V_0^n$,
- $\phi_n \circ f_n = \tilde{f}_n \circ \psi_n$ on $V_0^n - \cup_i D_i$.

Define $\phi_n := \psi_n$ outside $V_0^n$. The sequence of maps $\phi_n$ is eventually locally constant outside $K(f)$. So we obtained a limit map $\phi \colon \mathbb{C} \setminus K(f) \to \mathbb{C} \setminus K(\tilde{f})$, with $\phi(\partial V_0^0) = \partial \tilde{V}_0^0$, which is a $K$-quasiconformal conjugacy. By proposition 6.2, $\phi$ have a continuous extension to $K(f)$. □

*Proof of Theorem 6.* By proposition 11.1, there exists a map $\phi \colon \mathbb{C} \to \mathbb{C}$ which is a conjugacy between $f$ and $\tilde{f}$, with $\phi(V_0^0) = \tilde{V}_0^0$. Moreover, $\phi$ is $O(\delta)$-quasiconformal outside $K(f)$. By proposition 10.1, the critical pieces in the principal nest of $f$ and $\tilde{f}$ have bounded geometry and satisfies the complex bounds property (see Theorem 1), so, by Theorem 1, $\phi$ is quasiconformal in the complex plane. Because $\tilde{f}$ does not support invariant line fields in $K(\tilde{f})$, $\overline{\partial}\phi = 0$ on $K(\tilde{f})$, so $\phi$ is $O(\delta)$-quasiconformal in the complex plane. □

## 12. Infinitesimal theory

One of the most elegant and useful tools in Lyubich [Lyu99] and Avila, Lyubich and de Melo [ALdM] is the theory of infinitesimal perturbations of analytic maps. In these section we give an introduction to this theory, including some minor simplifications and straightforward generalizations that are appropriate in our case.

Let $f$ be a generalized polynomial-like map and $\tilde{U}$ a domain compatible with $f$. The **horizontal subspace** $E_f^h$ are the vector in $T\mathcal{B}_{nor}(\tilde{U})$ ($E_f^h$ depends on $\tilde{U}$, but we will omit this in the notation) so that there exists a quasiconformal vector field $\alpha \colon \mathbb{C} \to \mathbb{C}$, with $\overline{\partial}\alpha = 0$ on $K(f)$, with

$$v = \alpha \circ f - Df \cdot \alpha$$

in a neighborhood of $K(f)$.

The **vertical subspace** $E_f^v$ are the vectors $v \in T\mathcal{B}_{nor}(\tilde{U})$ so that there exists a holomorphic vector field $\alpha$ defined in $\overline{\mathbb{C}} - K(f)$, with $\alpha(\infty) = 0$, satisfying

(7) $$v = \alpha \circ f - Df \cdot \alpha$$

in some neighborhood of $K(f)$ in $\overline{\mathbb{C}} - K(f)$. The horizontal and vertical spaces were introduced by M. Lyubich [Lyu97]. If $|\cdot|_{\overline{\mathbb{C}}}$ denotes the Riemanian spherical metric and $\beta$ is a vector field on $\overline{U} \subset \overline{\mathbb{C}}$, define $|\beta|_{sph\ \overline{U}} := sup\{|\beta(x)|_{\overline{\mathbb{C}}} \colon x \in \overline{U}\}$.

**Proposition 12.1** ([Lyu99]). *Let $f \in \mathcal{B}_{nor}(U)$ so that $f$ has a generalized polynomial-like extension compatible with $U$. Let $V \Subset U$ be a neighborhood of $K(f)$ so that* $f^{-1}(V) \Subset V$. Then there exists $C$ so that, for any function $g$ close to $f$ so that $0$ does not escape of $V$, and $v = \alpha \circ g - Dg \cdot \alpha \in E_g^v$, we have

$$|\alpha|_{sph\ \overline{\mathbb{C}}\setminus V} \leq C|v|_{\mathcal{B}(V)}.$$

**Corollary 12.1.** *The vertical space $E_f^v$ is finite dimensional in any Banach slice compatible with $f$. Furthermore $\dim E_f^v$ does not depend on the Banach slice.*

*Proof.* It is sufficient to prove that the unity ball on $E_f^v$ is compact. Indeed, let $v_n = \alpha_n \circ f - Df \cdot \alpha_n \in E_f^v$, with $|v|_{\mathcal{B}(U)} \leq 1$. Let $V$ be as in proposition 12.1 and define $V_i := f^{-i}V$. Then there exists a sequence $C_i$ so that

$$|\alpha_n|_{sph\ \overline{\mathbb{C}}\setminus V_i} \leq C_i$$

In particular, there is a subsequence $n_k$ and a conformal vector field $\alpha_\infty$ on $\overline{\mathbb{C}}\setminus K(f)$ so that $\alpha_i \to \alpha_\infty$, where the convergence is uniform on compact sets in $\overline{\mathbb{C}} \setminus K(f)$. Furthermore, we can assume that $v_{n_k}$ converges uniformly to a conformal vector field $v_\infty$ in a neighborhood of $K(f)$ (for instance, in $V_1$). So we have

$$v_\infty = \alpha_\infty \circ f - Df \cdot \alpha_\infty$$

on $V_1 \setminus K(f)$. Consider $\tilde{U}$, $U \Subset \tilde{U}$, so that $f$ has an analytical extension to $\tilde{U}$. Use eq. 7 to extend $v_i$ and $v_\infty$ to $\tilde{U}$. By the maximum principle and proposition 12.1,

$$|v_{n_k} - v_\infty|_{\mathcal{B}(\tilde{U})} \leq C_1|\alpha_{n_k} - \alpha_\infty|_{sph\ \overline{\mathbb{C}}\setminus \tilde{U}} \leq C_1|\alpha_{n_k} - \alpha_\infty|_{sph\ \overline{\mathbb{C}}\setminus V_1} \leq C_2|v_{n_k} - v_\infty|_{\mathcal{B}(V_1)},$$

which finish the proof: the invariance of the dimension is an immediate consequence of equation 7, since any vertical vector defined in a domain $U$ has an analytical extension to a larger domain $\tilde{U}$. □

Let $f$ be a generalized polynomial-like map and let $v$ be a domain compatible with $f$. Denote $P(f) := \overline{\{f^i(0) \colon i > 0\}}$. If $v \in T\mathcal{B}(U)$, we will say that a quasiconformal vector fields $\alpha$ is $(f, v)$**-equivariant** on $P(f)$ if $v(0) = \alpha(f(0))$ and

$$\alpha(x_0) = lim_{x \to x_0} \frac{v(x) - \alpha(f(x))}{Df(x)},$$

for $x_0 \in P(f)$. The concept of equivariance was introduced in Avila, Lyubich & de Melo [ALdM]. The right side of the equation is well defined at the critical point (the unique problematic point) because $v(0) = \alpha(f(0))$ and $\alpha$ is quasiconformal: see [ALdM]. Indeed, always we have $lim_{x \to 0} (v(x) - \alpha(f(x)))/Df(x) = 0$. In particular, if $0 \in P(f)$ and $\alpha$ is $(f, v)$-equivariant on $P(f)$, then $\alpha(0) = 0$.

The Sullivan pullback argument has an infinitesimal counterpart:

**Proposition 12.2** (infinitesimal pullback argument[ALdM]). *Let $f$ be a generalized polynomial like map compatible with the domain $U$, whose critical point does not escape. Let $V \Subset U$ be so that $f \colon f^{-1}V \to V$ is a representation of $f$. Then there exists $\epsilon > 0$ and $C > 0$ so that, if $|g - f|_{\mathcal{B}(U)} \leq \epsilon$ and $g \colon g^{-1}V \to V$ is generalized polynomial like maps whose critical point also does not escape, then, for any vector $v \in T\mathcal{B}_{nor}$ that admits an $(g, v)$-equivariant quasiconformal vector field in $P(f)$, there exists a quasiconformal vector field $\alpha \colon \overline{\mathbb{C}} \to \overline{\mathbb{C}}$, $\alpha(\infty) = \alpha(1) = 0$, so that $v = \alpha \circ g - Dg \cdot \alpha$ on $\overline{V}$ and*

$$|\overline{\partial}\alpha|_{V\setminus K(g)} \leq C|v|_{\mathcal{B}(U)}.$$

**Corollary 12.2** (semicontinuity of the horizontal subspaces[ALdM])**.** *Let $f$, $f_n$ be generalized polynomial-like maps which are compatible with the Banach slice $\mathcal{B}_{nor}(U)$ and so that $f_n \to f$ in $\mathcal{B}_{nor}(U)$. Assume that $f$ does not support invariant line fields in $K(g)$. If $v_n \in E_{f_n}^h \subset T\mathcal{B}(U)$ and $v_n \to v$ in $T\mathcal{B}_{nor}(U)$, then $v \in E_f^h$. In particular $E_f^h$ is closed.*

The infinitesimal version of the Douady & Hubbard decomposition of the dynamics of a polynomial-like map in the hybrid class and its external map is the following:

**Proposition 12.3** (Lyubich Decomposition[Lyu99])**.** *In a Banach slice $\mathcal{B}_{nor}(U)$ compatible with $f$, $T\mathcal{B}_{nor}(U) = E_f^h \oplus E_f^v$.*

M. Lyubich stated the above decomposition for quadratic polynomial-like maps (connected domain), but the proof can be carry out as well in our case.

The proof of corollary 12.1 does not tell us anything about the finite dimension of $E_f^v$. The following result will be very useful to solve this problem. The proof of it is a variation of a result in [Lyu99].

**Proposition 12.4.** *Let $f$ be a generalized polynomial, with non escaping critical points. Assume there are not invariant line fields on $K(f)$ and let $U$ be a domain compatible with $f$. Then there exists $\epsilon > 0$ so that, for any generalized polynomial like map $g$, with non escaping critical points and $\epsilon$-close to $f$ on $\mathcal{B}_{nor}(U)$, we have $\dim E_f^v = \dim E_g^v$.*

*Proof.* Firstly, suppose that there exists $f_n \to f$ on $\mathcal{B}_{nor}(\tilde{U})$, with $\dim E_f^v > \dim E_{f_n}^v$. Then, for every $n$ we can find $v_n \in T\mathcal{B}_{nor}(U)$ so that $v_n \in E_f^v \cap E_{f_n}^h$ and $|v_n|_{\mathcal{B}(U)} = 1$. Since $E_f^v$ is finite dimensional, we can assume that $v_n \to v \in E_f^v$, $v \neq 0$. But, by the semi-continuity of the horizontal spaces, $v \in E_f^h$, what is a contradiction. To finish the proof, it is sufficient so prove that, if there exists a sequence of maps $f_n \to f$ so that $\dim E_{f_n}^v \geq N$, then $\dim E_f^v \geq N$. Indeed, by a classical lemma by Riesz, there exist "almost orthonormal" families of vectors $\{v_n^1, \ldots, v_n^N\} \subset E_{f_n}^v$. In other words, $|v_n^i|_{\mathcal{B}(U)} = 1$ and, if we denote

$$span\{u_1, \ldots, u_N\} := \{\sum c_i u_i \colon c_i \in \mathbb{C}\},$$

we have

$$dist(v_n^i, span\{v_n^1, \ldots, v_n^{i-1}\})_{\mathcal{B}(U)} \geq \frac{1}{2},$$

for $1 < i \leq N$. Let $\hat{U} \Subset U$. Since $|v_n^i|_{\mathcal{B}(U)} = 1$, we can assume that $v_\infty^i \to v_\infty^i$ on $T\mathcal{B}_{nor}(\hat{u})$, for some vectors $v_n^\infty$. By proposition 12.1,

- $v_\infty^i \in E_f^v \setminus \{0\}$ and
- $dist(v_\infty^i, span\{v_\infty^1, \ldots, v_\infty^{i-1}\})_{\mathcal{B}(U)} \geq C > 0$, for $1 < i \leq N$.

So $\dim E_f^v \geq N$. □

Denote $\mathcal{O}(x) := \{f^i(x) \colon i \geq 0\}$. The following is a slight variation of a similar result to maps with strictly preperiodic critical point in [ALdM]:

**Proposition 12.5.** *Let $f$ be a generalized polynomial-like map with an unique critical point, which is $N$-periodic. Then, if $U$ is compatible with $f$, then $v \in E_f^h$ if and only if there is a vector $\alpha \colon \mathcal{O}(0) \to \mathbb{C}$ so that $v = \alpha \circ f - Df \cdot \alpha$, with $\alpha(0) = 0$. In particular, codim $E_f^h = 1$.*

*Proof.* It is obvious that the assumptions under $\alpha$ are necessary. We will prove that they are sufficient. Replace $\alpha$ by a quasiconformal vector field on $\overline{\mathbb{C}}$ that satisfies the same hypothesis. Then $\alpha$ is equivariant with respect to $(f, v)$ on $P(f)$. We can also assume that $\alpha$ is analytic in a neighborhood of $\mathcal{O}(0)$. It is easy to see in the proof of the infinitesimal pullback argument[ALdM] that there exists a quasiconformal vector field $\alpha\colon \mathbb{C} \to \mathbb{C}$, analytic in the interior of $K(f)$, so that $v = \alpha \circ f - Df \cdot \alpha$ on a neighborhood of $K(f)$. Since $m(J(f)) = 0$, $v \in E_f^h$. To prove that $codim\ E_f^h = 1$, denote $v_i = v(f^i(0))$ and $\alpha_i = \alpha(f^i(0))$. Then

$$v_i = \alpha_{i+1\ mod\ n} - Df(f^i(0)) \cdot \alpha_i.$$

Since the determinant of this system of equations is 1, for each $N$-uple $(v_0, \ldots, v_{N-1})$ there exists an unique solution $\alpha_j(v_0, \ldots, v_{N-1})$. By the above characterization of $E_f^h$, the horizontal space is the kernel of the non zero functional

$$v \to \alpha_0(v(0), \ldots, v(f^{N-1}(0))).$$

$\square$

**Corollary 12.3.** *Let $f$ be a generalized polynomial-like map with an unique critical point, assume that there are not invariant line fields on $K(f)$ and let $U$ be a domain compatible with $f$. Assume that there exists a sequence if generalized polynomial-like maps with representations in $\mathcal{B}_{nor}(U)$ so that $f_n$ has a superattractor and $f_n \to f$ on $\mathcal{B}_{nor}(U)$. Then $\dim E_f^v = 1$.*

*Proof.* Propositions 12.4 and 12.5. $\square$

A generalized polynomial-like map $f_1\colon \cup_i U_1 \to V_1$ is **hybrid conjugated** to another $f_2\colon \cup_i U_2 \to V_2$ if there exists a quasiconformal map $h\colon \mathbb{C} \to \mathbb{C}$ so that

- $h(V_1) = V_2$,
- $h \circ f_1 = f_2 \circ h$ on $V_1$,
- $\overline{\partial} h$ vanishes at $K(f_1)$.

The **hybrid class** of $f_1$ is the set of generalized polynomial-like maps which are hybrid conjugated with $f_1$.

**Corollary 12.4.** *Suppose that $f_1$ has an unique critical point, which does not escape, and $K(f_1)$ does not support invariant line fields. Assume that $f_2$ is a generalized polynomial-like map which are hybrid conjugated with $f_2$. Then the codimension of $E_{f_2}^h$ in an arbitrary Banach slice compatible with $f_2$ is equal to the codimension of $E_{f_1}^h$.*

*Proof.* Indeed, let $h$ be a quasiconformal homeomorphism between to generalized polynomial-like maps $f_a$ and $f_b$ and the hybrid class. So $h$ induce a Beltrami path $h_\lambda \circ f_a \circ h_\lambda^{-1}$ ($h_0 = Id$ and $h_1 = h$). Given two Banach slices $\mathcal{B}(U_a)$ and $\mathcal{B}(U_b)$ compatible with $f_a$ and $f_b$, we can find, using the compactness of $[0, 1]$ and proposition 12.4, a sequence of Banach slices $U_a = \hat{U}_0, \ldots, \hat{U}_n = U_b$ and closed intervals $I_i \subset [0, 1]$, with $0 \in I_0$ and $1 \in I_n$, $I_i \cap I_{i+1} \neq \phi$, so that

- For $\lambda \in I_i$, $f_\lambda$ is compatible with $\mathcal{B}_{nor}(\hat{U}_{i-1})$, $\mathcal{B}_{nor}(\hat{U}_i)$ and $\mathcal{B}_{nor}(\hat{U}_{i+1})$.
- $codim\ E_{f_{\lambda_0}}^h = codim\ E_{f_{\lambda_1}}^h$ on $\mathcal{B}_{nor}(\hat{U}_i)$, for $\lambda_0, \lambda_1 \in I_i$.

which completes the proof. $\square$

## 13. Induced maps

In this section we study induced transformations in a quite general setting: induced transformations includes the classical renormalization operator and the Fibonacci renormalization operator.

Let $f$ be a generalized polynomial-like map, $\mathcal{B}(U)$ be a Banach slice compatible with $f$ and let $\beta_g$ be a non zero analytical function defined in a neighborhood of $f$ in $\mathcal{B}(U)$ ($\beta_g$ will be called **normalization function** or simply 'normalization'). Consider $U_1, \ldots, U_n$ be simple connected domains and $j_1, \ldots, j_n \in \mathbb{N}$ so that

- the closures of $U_i$ are disjoint.
- the sets $\overline{U_i}, \ldots, f^{j_i}(\overline{U_i})$ are contained in $U$.

Define $\hat{U}_i = 1/\beta_f U_i$. Then we define a map of a neighborhood of $f$ in $\mathcal{B}(U)$ in the Banach slice $\mathcal{B}(\hat{U}_1 \cup \cdots \cup \hat{U}_n)$ in the following way: for each $g$ in a small neighborhood of $f$ in $\mathcal{B}(U)$ define $\mathcal{R}g \in \mathcal{B}(\hat{U}_1 \cup \cdots \cup \hat{U}_n)$ as

$$\mathcal{R}g(x) = \frac{1}{\beta_g} g^{j_i}(\beta_g x), \text{ for } x \in \hat{U}_i.$$

The map $\mathcal{R}$ is called the **(normalized) induced transformation** associated with the $2n$-uple $(U_1, j_1, \ldots, U_n, j_n)$ and the normalization $\beta$. The set $U_1 \cup \cdots \cup U_n$ is the **support** of the induced transformation.

Let $v \in \mathcal{B}(U)$ be a vector. Denote by $a_n$ the following functions:

$$a_1(x) := v(x),$$
$$a_n(x) := v(f^{n-1}(x)) + Df(f^{n-1}(x))a_{n-1}(x).$$

Then (we replaced $\beta_f$ by $\beta$ to simplify the notation)

$$(D\mathcal{R}_f \cdot v)(x) = \frac{\partial}{\partial t}\bigg|_{t=0} \frac{1}{\beta_{f+tv}}(f+tv)^n(\beta_{f+tv}x)$$
$$= \frac{-D\beta \cdot v}{\beta}\frac{1}{\beta}f^n(\beta x) + \frac{1}{\beta}\big(a_n(\beta x) + Df^n(\beta x)(D\beta \cdot v)x\big).$$

This is an ugly expression, but in the horizontal space we have a much nicer description of the action of $D\mathcal{R}$:

**Proposition 13.1.** *Let $\mathcal{R}$ be an normalized induced transformation with normalization $\beta_f$. If $v = \alpha \circ f - Df \cdot \alpha \in E_f^h$ then*

$$D\mathcal{R}_f \cdot v = r(\alpha) \circ \mathcal{R}f - D(\mathcal{R}f) \cdot r(\alpha),$$

*where*

$$r(\alpha) := \frac{1}{\beta_f}\alpha(\beta_f x) - \frac{1}{\beta_f}(D\beta_f \cdot v) \cdot x.$$

*Furthermore, if $\beta_f$ is the analytic continuation of a non parabolic period point, then*

$$r(\alpha) := \frac{1}{\beta_f}\alpha(\beta_f x) - \frac{1}{\beta_f}\alpha(\beta_f) \cdot x.$$

*Proof.* So assume $v = \alpha \circ f - Df \cdot \alpha$. Replacing this expression in the above equation we obtain, by induction, that, for $n > 1$:

$$a_n = \alpha \circ f^n - Df^n \cdot \alpha.$$

It follow that $D\mathcal{R} \cdot v$ can be rewritten as $D\mathcal{R}_f \cdot v = r(\alpha) \circ \mathcal{R}f - D(\mathcal{R}f) \cdot r(\alpha)$, where

$$r(\alpha)(x) := \frac{1}{\beta}\alpha(\beta x) - \frac{1}{\beta}(D\beta \cdot v)x.$$

If we choice $\beta_f$ as the as the analytic continuation of a non parabolic periodic point of $f_0$ with period $i$, then we have $f^i(\beta_f) = \beta_f$, so

$$D\beta_f \cdot v = a_i(\beta) + Df^i(\beta)D\beta_f \cdot v = \alpha(\beta) - Df^i(\beta)\alpha(\beta) + Df^i(\beta)D\beta_f \cdot v,$$

hence $D\beta_f \cdot v = \alpha(\beta)$ and

$$r(\alpha)(x) := \frac{1}{\beta}\alpha(\beta x) - \frac{1}{\beta}\alpha(\beta)x.$$

In particular, $v(1) = r(1) = 0$. $\square$

All items in the proposition below, except the last one, are obvious generalizations of similar statements to the classic renormalization operator in de Melo & van Strien book[dMvS], Lyubich[Lyu99] and Avila, Lyubich & de Melo [ALdM] (in remark 7.2 in [ALdM], Avila, Lyubich & de Melo stated item 2). Item 3 is consequence of proposition 13.1. The statement about independence of induced transformations will be useful only when we are dealing with generalized polynomial-like maps with many critical points, and it will be not used here. Its proof is a straightforward generalization of the proof of the density of the image of the derivative of the classic renormalization operator [ALdM].

**Proposition 13.2.** *Let $\mathcal{R}$ be a normalized induced transformation associated with the $2n$-uple $(U_1, j_1, \ldots, U_n, j_n)$, defined in a neighborhood of $f$ in a compatible Banach slice $\mathcal{B}_{nor}(U)$ in $\mathcal{B}_{nor}(U_1 \cup \cdots \cup U_n)$. Then*

1. **Compactness.** *The induced transformation is a compact transformation.*
2. **Derivative with dense image.** *The derivative $D\mathcal{R}$ is a compact linear transformation with dense image.*
3. **Invariance of horizontal subspaces.** *Assume that $f$ and $\mathcal{R}f$ have generalized polynomial-like extensions so that $\mathcal{B}_{nor}(U)$ and $\mathcal{B}_{nor}(U_1 \cup \cdots \cup U_n)$ are compatible with them. Then $D\mathcal{R}_f(E_f^h) \subset E_{\tilde{g}}^h$.*
4. **Non-singularity in the vertical directions.** *Under the conditions of the previous item, assume that $D\mathcal{R}^{-1}E_g^h \subset E_f^h$: Then*

$$\tilde{D\mathcal{R}}_f \colon T\mathcal{B}_{nor}(U)/E_f^h \to T\mathcal{B}_{nor}(U_1 \cup \cdots \cup U_n)/E^h\tilde{g}$$

   *is invertible.*

*Moreover we have*

- **Independence between induced transformations.** *If $\mathcal{R}_i$, $i \leq i_0$, are induced transformations associated with the $2n$-uples*

$$(U_1^i, j_1^i, \ldots, U_n^i, j_{n_i}^i),$$

   *where $U_j^i \subset U$, and $\mathcal{R}_i$ have two by two disjoint supports, then*

$$\{(D\mathcal{R}_1(f) \cdot v, \ldots, D\mathcal{R}_i(f) \cdot v, \ldots, D\mathcal{R}_{i_0}(f) \cdot v) \colon v \in \mathcal{B}(U)\}$$

   *is dense in $T\mathcal{B}_{nor}(U_1^1 \cup \cdots \cup U_{n_1}^1) \oplus \cdots \oplus T\mathcal{B}_{nor}(U_1^i \cup \cdots \cup U_{n_i}^i) \oplus \cdots \oplus T\mathcal{B}_{nor}(U_1^{i_0} \cup \cdots \cup U_{n_{i_0}}^{i_0}).*

## 14. Contraction on the hybrid class

Let $\tilde{f}_i \colon \tilde{V}_0^1 \cup \tilde{V}_1^1 \to \tilde{V}_0^0$, $i = 1, 2$, be the generalized polynomial-like extensions of $\tilde{f}_i$, as previously defined. Denote by $\mathcal{T}_{Fib}$, called **Fibonacci tower** [vSN94], the indexed family of generalized polynomial-like maps $\mathcal{R}^i \tilde{f}_1 \colon \beta^i \tilde{V}_0^1 \cup \beta^i \tilde{V}_1^1 \to \beta^i \tilde{V}_0^0$, with $i \in \mathbb{Z}$, where either $\mathcal{R}^i \tilde{f}_1(z) := 1/\beta^i \tilde{f}_1(\beta^i \cdot z)$, if $i$ is even, or $\mathcal{R}^i \tilde{f}_1(z) := 1/\beta^i \tilde{f}_2(\beta^i \cdot z)$, if $i$ is odd. This notation is coherent with the previous definitions of the Fibonacci operator $\mathcal{R}$, since this operator is injective. Let $\mu$ be a Beltrami field in $\mathbb{C}$. We say that $\mu$ is invariant by the tower $\mathcal{T}$ if

$$(\mathcal{R}^i \tilde{f}_1)_* \mu = \mu,$$

for $i \in \mathbb{Z}$. We say that the Fibonacci tower $\mathcal{T}_{Fib}$ does not support invariant line fields if every Beltrami field which is invariant by $\mathcal{T}_{Fib}$ vanishes almost everywhere in $\mathbb{C}$.

**Proposition 14.1** ([vSN94]). *The Fibonacci tower $\mathcal{T}_{Fib}$ does not support invariant line fields.*

The following result can be proved using arguments as in van Strien & Nowicki [vSN94]:

**Proposition 14.2** (Topological convergence[vSN94]). *Let $f \colon V_0^1 \cup V_1^1 \to V_0^0$ be a generalized polynomial-like map which admits a $K$-quasiconformal map $h \colon \mathbb{C} \to \mathbb{C}$ which is a conjugacy between $f$ and $\tilde{f}$, with $h(V_0^0) = \tilde{V}_0^0$. Then there exists $N = N(K)$ so that the restriction $f_n \colon V_0^{n+1} \cup V_1^{n+1} \to V_0^n$ of the first return map to $V_0^n$ (recall the notation in section 4.3) satisfies the following: if $\tau_n := h^{-1}(\beta^n)$, then*

$$\frac{1}{\tau_n} f_n(\tau_n \cdot z)$$

*has an analytical extension to $\mathcal{B}_{nor}(U)$, provided that $n \geq N(f)$, and*

$$|\tau_n^{-1} f_n(\tau_n \cdot z) - \tilde{f}|_{\mathcal{B}(U)} \to 0.$$

The main result in this section is that this convergence is, indeed, exponential. Let $\mathcal{A}(K, V, \epsilon)$ be the set of generalized polynomial like maps $f$ in $\mathcal{B}_{nor}(U)$ so that

- $\mathcal{R}^n f$ is defined for $n \geq 0$ and

$$|\mathcal{R}^n f - \tilde{f}|_{\mathcal{B}(U)} \leq \epsilon.$$

- There exists a $K$-quasiconformal map $\phi \colon \mathbb{C} \to \mathbb{C}$ so that $\phi(V) \subset \overline{U}$ and $\phi \circ \tilde{f} = f \circ \phi$ on $V$.

Note that the (closed) set $\mathcal{A}(K, V_0^0, \epsilon)$ is invariant by the action of $\mathcal{R}$, if $\epsilon$ is small enough. Moreover $\overline{\mathcal{R}\mathcal{A}(K, V_0^0, \epsilon)}$ is a compact set.

**Corollary 14.1.** *Given $\epsilon$ small enough, there exists $K_\epsilon = O(\epsilon)$ so that if $f \colon V_0^1 \cup V_1^1 \to V_0^0$ is a generalized polynomial-like map which admits a $\tilde{K}$-quasiconformal map $h \colon \mathbb{C} \to \mathbb{C}$ which is a conjugacy between $f$ and $\tilde{f}$, with $h(V_0^0) = \tilde{V}_0^0$, then, for $n \geq N(\tilde{K})$, $\mathcal{R}^n f \in \mathcal{B}_{nor}(U)$ and*

$$\mathcal{R}^n f \in \mathcal{A}(K_\epsilon, V_0^0, \epsilon).$$

*Proof.* Given $\epsilon$, by Theorem 6, there exists $K = O(\epsilon)$ so that if $f$ satisfies

$$|\mathcal{R}^n f - \tilde{f}| \leq \epsilon$$

for $n \geq 0$, then there exists hybrid conjugacy between $f\colon V_0^1 \cup V_1^1 \to V_0^0$ and $\tilde{f}\colon \tilde{V}_0^1 \cup \tilde{V}_1^1 \to \tilde{V}_0^0$ which is $K$-quasiconformal in the complex plane. Proposition 14.2 complete the proof. $\square$

**Remark 14.1.** *The proof of a result similar to Corollary 14.1 is easier in the classic renormalization theory for quadratic polynomials: we do not know if it is possible to prove it without proposition 11.1. The problem is that there is not a canonical straightening for Fibonacci generalized polynomial-like maps.*

**Theorem 7** (Contraction in the horizontal direction). *Given $\mathcal{A}(K, V_0^0, \epsilon)$, with $\epsilon$ small enough, there exists $\lambda < 1$ and $N \in \mathbb{N}$ so that if $v \in E_f^h$ and $f \in \mathcal{A}(K, V, \epsilon)$, then*
$$|D\mathcal{R}_f^N \cdot v|_{\mathcal{B}(U)} \leq \lambda |v|_{\mathcal{B}(U)}.$$

*Proof.* Consider $f \in \mathcal{A} := \mathcal{A}(K, V_0^0, \epsilon)$ and $v \in E_f^h$, $|v|_{\mathcal{B}(U)} \leq 1$. Then, by the infinitesimal pullback argument, there exists a quasiconformal vector field $\alpha\colon \overline{\mathbb{C}} \to \overline{\mathbb{C}}$ so that $v = \alpha \circ f - Df \cdot \alpha$ on $\overline{V}_0^1 \cup \overline{V}_1^1$, with $|\overline{\partial}\alpha|_\infty \leq C$. This implies, by the definition of the Fibonacci renormalization operator, proposition 12.2 and proposition 13.1, that, for $n \geq 1$,

(8) $$D\mathcal{R}_f^n \cdot v = \alpha_n \circ \mathcal{R}^n f - D(\mathcal{R}^n f) \cdot \alpha_n,$$

on $U = U_1 \cup U_2$, where

$$\alpha_n(z) := \frac{1}{\beta_{\mathcal{R}^{n-1}f} \cdots \beta_f} \alpha(\beta_{\mathcal{R}^{n-1}f} \cdots \beta_f \cdot z) - \frac{1}{\beta_{\mathcal{R}^{n-1}f} \cdots \beta_f} \alpha(\beta_{\mathcal{R}^{n-1}f} \cdots \beta_f) \cdot z.$$

Note that $|\overline{\partial}\alpha_n|_\infty = |\overline{\partial}\alpha|_\infty \leq C$ on $\overline{\mathbb{C}}$, for $n \geq 1$. Recall that $C$ does not depend on $f$ or $v$. Since the set of $C$-quasiconformal vectors which vanishes in three points (in our case, 0, 1 and $\infty$) in the Riemann sphere is a compact set relative to the sup norm induced by the spherical metric in the Riemann sphere, There exists $C > 0$ so that
$$|D\mathcal{R}_f^n \cdot v|_{\mathcal{B}(U)} \leq C$$
for every $f \in \mathcal{A}(K, V, \epsilon)$, $v \in E_f^h$, $|v| \leq 1$. Since $\mathcal{R}$ is a compact operator, we have that
$$S := \overline{\{(\mathcal{R}f, D\mathcal{R}_f \cdot v)\colon f \in \mathcal{A}, v \in E_f^h, |v|_{\mathcal{B}(U)} \leq 1\}}$$
is a compact set. Furthermore, since the horizontal subspaces are semi-continuous (corollary 12.2), we get $S \subset \{(f, v)\colon f \in \mathcal{A}, v \in E_f^h\}$.

We claim that, for any $v \in E_f^h$, $f \in \mathcal{A}$, we have
$$|D\mathcal{R}_f^n \cdot v|_{\mathcal{B}(U)} \to 0.$$
It is sufficient to prove that $\alpha_n \to 0$. Take any convergent subsequence $\alpha_{n_k} \to \alpha_\infty$. In particular, $\overline{\partial}\alpha_{n_k}$ converges to $\overline{\partial}\alpha_\infty$ in the distributional sense. Denote $\beta_i := \beta_{\mathcal{R}^i f}$. Note that the Beltrami field $\overline{\partial}\alpha_n$ is invariant by the maps

$$\mathcal{R}^n f, \; \frac{1}{\beta_{n-1}} \mathcal{R}^{n-1} f(\beta_{n-1} \cdot), \ldots, \frac{1}{\beta_{n-1} \cdots \beta_0} f(\beta_{n-1} \cdots \beta_0 \cdot).$$

Since $\mathcal{R}^i f \to_i \tilde{f}$, this sequence of maps converges to the Fibonacci tower, which implies that $\overline{\partial}\alpha_\infty$ is invariant by the Fibonacci tower. Since the Fibonacci tower does not admit invariant line fields, $\alpha_\infty$ is a conformal vector field in the Riemann sphere. So $\alpha_\infty$ is zero (because $\alpha_\infty$ is zero in three points in the Riemann sphere: 0, 1 and $\infty$). Fix $\delta$ small. Then, for every $(f,v) \in S$, there exists $N_{(f,v)}$ so that

$$|D\mathcal{R}_f^n \cdot v|_{\mathcal{B}(U)} \leq \delta,$$

for $n \geq N_{(f,v)}$. Since $|D\mathcal{R}_f^n| \leq C$, we can find a neighborhood of $B_{(f,v)}$ of $(f,v)$ in $S$ so that

$$|D\mathcal{R}_g^n \cdot u|_{\mathcal{B}_{nor}(U)} \leq 2C\delta,$$

for $(g,u) \in B_{(f,v)}$ and $n \geq N_{(f,v)}$. Since $S$ is compact, this finish the proof. □

**Corollary 14.2** (Contraction on the hybrid class). *There exists $\lambda < 1$ so that, if $f$ is a generalized polynomial-like map which is quasiconformally conjugated with one of the maps in the Fibonacci cycle, then, for $n > N = N(f)$, the nth generalized renormalization of $f$ have a representation in $\mathcal{B}_{nor}(U)$ and furthermore*

$$|\mathcal{R}^n f - \tilde{f}|_{\mathcal{B}(U)} \leq \lambda^n.$$

*Proof.* Let $\phi \colon \mathbb{C} \to \mathbb{C}$ be a quasiconformal map which is a conjugacy between $f$ and $\tilde{f}$ in a neighborhood of $J(\tilde{f})$. Consider the Beltrami path induced by $\phi$: in other words, $f_t := \phi_t \circ \tilde{f} \circ \phi_t^{-1}$, $|t| \leq 1$, where

$$\frac{\overline{\partial}}{\partial}\phi_t = t \cdot \frac{\overline{\partial}}{\partial}\phi.$$

By the topological convergence of $\mathcal{R}^n f$ to $\tilde{f}$ (see [vSN94]), given $\epsilon$, there exists $N = N(f, \epsilon)$ so that $\mathcal{R}^n f_t \in \mathcal{B}_{nor}(U)$ and

$$|\mathcal{R}^n f_t - \tilde{f}|_{\mathcal{B}(U)} \leq \epsilon,$$

for $|t| \leq 1$. In particular, by proposition 11.1 and Theorem 1, there exists $K = K(\epsilon)$ and a definitive neighborhood $V$ on $J(\tilde{f})$ so that $\mathcal{R}^n f_t \in \mathcal{A}(K_\epsilon, V_0^0, \epsilon)$, for $n > N_f$ and $|t| \leq 1$. Since

$$\frac{d\mathcal{R}^n f_t}{dt} \in E_{\mathcal{R}^n f_t}^h,$$

we can use Theorem 7 to finish the proof. □

The above argument works also to the classic renormalization operator [Lyu99] (even for non periodic combinatorics), the multimodal renormalization operator [Sm01] , and possibly in any other situation in complex dynamics where the McMullen arguments with towers can be applied, as the renormalization theory for critical homeomorphisms or coverings of the circle.

## 15. Hyperbolicity

**Proposition 15.1.** *There exists a sequence of the generalized polynomial-like maps $f_{n,i}$, $i = 1, 2$, real in the real line, so that*

- *The map $f_{n,i}$ has a superattractor with period $S_n$ (the Fibonacci sequence),*
- *The sequence of closest returns of the critical point is exactly $f_{n,i}^{S_i}(0)$, $i < n$.*
- *The maps $f_{n,i}$ are compatible with $U$ and $f_{n,i} \to \tilde{f}_i$.*

*In particular, codim $E_{\tilde{f}}^h = 1$.*

*Proof.* Consider the Fibonacci parameter $c_\infty$ in the family $x^d + c$. Select a generalized polynomial-like extension $f_{n_0}^{c_\infty} \colon \hat{U}_1 \cup \hat{U}_2 \to \hat{V}$ of a generalized renormalization of $x^d + c_\infty$, which has the Fibonacci combinatorics. There exists an unique sequence of parameters $c_n$ so that

- $x^d + c_n$ has a critical point with period $S_n$ (the Fibonacci sequence),
- The sequence of the closest return times of the critical point is $S_i$, $i < n$.

Because any accumulation point of the sequence $c_n$ is a parameter with Fibonacci combinatorics, we have $c_n \to c_\infty$, since that, by the bounded geometry of the postcritical set, the non existence of invariant line fields at the Julia of Fibonacci polynomials and the Sullivan pullback argument, $c_\infty$ is the unique such parameter (see [vSN94]). It is easy to see that $P(x^d + c_n) \to P(x^d + c_\infty)$. In particular, the $n_0$-th generalized renormalization of $x^d + c_n$), denoted $f_{n_0}^{c_n}$, has also a polynomial-like extension with, say, the same image $\hat{V}$ and domains $\hat{U}_0^n$ and $\hat{U}_1^n$, where these domains are very close to $\hat{U}_0$ and $\hat{U}_1$. This is sufficient to conclude that *codim* $E_{\tilde{f}}^h = 1$, since $f_{n_0}$ is in the hybrid class of $\tilde{f}_1$ (or $\tilde{f}_2$). Now we could use the dynamics of the Fibonacci renormalization to conclude the proof, but we will use a more general argument: just to fix the notation, suppose that $f_{n_0}$ is hybrid conjugated with $\tilde{f}_1$. Then select a generalized polynomial-like extension $\tilde{f}_1 \colon \tilde{U}_0 \cup \tilde{U}_1 \to \tilde{V}$ which has a representation in $U$. Then there exists $K$-quasiconformal maps, symmetric with respect the real line,

$$\phi_n \colon \mathbb{C} - \hat{U}_0^n \cup \hat{U}_1^n \to \mathbb{C} - \tilde{U}_0 \cup \tilde{U}_1$$

so that

- $\phi_n(\partial \hat{U}_i^n) = \partial \tilde{U}_i^n$,
- $\phi_n \circ f_{n_0}^{c_n} = \tilde{f} \circ \phi_n$ on $\partial \hat{U}_0 \cup \partial \hat{U}_1$.

Let $\mu_n$ be the Beltrami field which coincides with $\overline{\partial} \phi_n$ on $\mathbb{C} - \hat{U}_0^n \cup \hat{U}_1^n$, it is invariant by $f_{n_0}^{c_n}$ and $\mu_n = 0$ on $K(f_{n_0}^{c_n})$. Let $\psi_n$ be the $K$-quasiconformal maps, symmetric with respect the real line, so that $\overline{\partial} \psi_n = \mu_n$. We claim that the generalized polynomial-like maps $\hat{f}_n := \psi_n \circ f_{n_0}^{c_n} \circ \psi_n^{-1}$ converges to $\tilde{f}_1$. Clearly any convergent subsequence converges to a generalized polynomial-like map with the form $\hat{f} := \psi \circ f_{n_0}^{c_\infty} \circ \psi^{-1}$. Note that the map $\phi_n \circ \psi_n^{-1} \colon \overline{\mathbb{C}} - \psi_n(\hat{U}_0^n \cup \hat{U}_1^n) \to \overline{\mathbb{C}} - (\tilde{U}_0 \cap \tilde{U}_1)$ is conformal. Moreover $\phi_n \circ \psi_n^{-1} \circ \hat{f}_n = \tilde{f}_1 \circ \phi_n \circ \psi_n^{-1}$ on $\psi_n(\partial \hat{U}_0 \cup \partial \hat{U}_1)$. Taking the limit we conclude that there exists $\phi \circ \psi^{-1}$ so that

- $\phi \circ \psi^{-1}$ is conformal outside $\psi(\hat{U}_0^n \cup \hat{U}_1^n)$,
- $\phi \circ \psi^{-1} \circ \hat{f} = \tilde{f}_1 \circ \phi \circ \psi^{-1}$ on $\psi(\partial \hat{U}_0 \cup \partial \hat{U}_1)$.

Construct a real symmetric quasiconformal map which coincides with a quasisymmetric conjugacy between $\hat{f}$ and $\tilde{f}_1$, in the real line, and with $\phi \circ \psi^{-1}$ on $\overline{\mathbb{C}} - \psi(\partial \hat{U}_0 \cup \partial \hat{U}_1)$. By the Sullivan pullback argument, we obtain a quasiconformal map $h \colon \overline{\mathbb{C}} \to \mathbb{C}$ which is conformal outside $K(\hat{f})$ and it is conjugacy between $\hat{f}$ and $\tilde{f}_1$. Since there are not invariant line field in $K(\hat{f})$, $h$ is an affine map, which finish the proof. □

**Proposition 15.2** (No small orbits). *There exists $\epsilon > 0$ so that if $f \in \mathcal{B}_{nor}(U)$ satisfies*

$$|\mathcal{R}^n f - \tilde{f}|_{\mathcal{B}(U)} \leq \epsilon$$

*for $n \geq 0$ then $f$ has a generalized polynomial like restriction which is hybrid conjugated with either $\tilde{f}_1$ or $\tilde{f}_2$ and*

$$|\mathcal{R}^n f - \tilde{f}|_{\mathcal{B}(U)} \leq C\lambda^n,$$

*for some $C > 0$ and $\lambda < 1$.*

*Proof.* Obvious consequence of Theorem 6 and Proposition 14.2. □

**Theorem 8.** *The Fibonacci cycle is hyperbolic with one dimensional unstable manifold.*

*Proof.* Note that, by proposition 12.3 and 15.1, *codim* $E^h_{\tilde{f}} = 1$. We already proved that $D\mathcal{R}^2_{\tilde{f}}$ keeps invariant and it is a contraction on the Banach subspace $E^h_{\tilde{f}}$. Note that the spectrum of $D\mathcal{R}^2_{\tilde{f}}$ cannot be contained inside the open disc, otherwise $\tilde{f}$ is an attracting fixed point, what implies, by proposition 11.1 and Theorem 1, that any function in $\mathcal{B}_{nor}(U)$ very close to $\tilde{f}$ admits a quasiconformal conjugacy with $\tilde{f}$, what is a contradiction with proposition 15.1. So, if $D\mathcal{R}^2_{\tilde{f}}$ is not hyperbolic, then the spectrum is contained in the closed unit disc and intersect the unit circle. By the Small Orbits Theorem [Lyu99], for any small $\epsilon$, there exists $f$ so that $|\mathcal{R}^i f - \tilde{f}| \leq \epsilon$, for $i \geq 0$, but $\mathcal{R}^i f$ does not converges exponentially to $\tilde{f}$, which is a contradiction with corollary 15.2. This finished the proof. □

## 16. Universality

In this section, we will restrict ourselves to the subspace of even maps (but we will keep the same notation to the Banach slices). This restriction is not really necessary in the universality result, but this assumption will simplify the exposition a little. Recall that the Fibonacci renormalization operator is defined in two balls $B_1$ and $B_2$ in the Banach slice $\mathcal{B}_{nor}(U_0 \cup U_1)$: The ball $B_i$ is a $\epsilon_0$-ball whose center is $\tilde{f}_i$. Consider the linear transformation $\phi\colon \mathcal{B}_{nor}(U_0 \cup U_1) \to \mathcal{B}_{nor}(U_0 \cup U_1)$ defined by $\phi(g) = \hat{g}$, where $\hat{g}$ coincides with $g$ on $\overline{U}_0$ and with $-g$ on $\overline{U}_1$. Then

- $\phi$ is a linear isometry of $\mathcal{B}_{nor}(U_0 \cup U_1)$,
- $\phi \circ \phi = Id$,
- $\phi \circ \mathcal{R} = \mathcal{R} \circ \phi$.

In particular $(\mathcal{R} \circ \phi)^2 = \mathcal{R}^2$. Since $\phi(\tilde{f}_i) = \tilde{f}_{1-i}$, we have that $\tilde{f}_i$ are hyperbolic fixed points of the operator $\mathcal{R} \circ \phi$. It is a question of taste to work with either $\mathcal{R}$ or $\mathcal{R} \circ \phi$: we prefer to use the first one in most of this work because its definition is more 'dynamical', what simplify the exposition. The following lemma can be proved with the argument applied in [ALdM] to prove that the set of real analytic maps which are topologically conjugated with a given map $f$ is a real analytic manifold (indeed, our situation is simpler, since we will restrict ourselves to maps with a superattractor):

**Proposition 16.1.** *There exist $\epsilon > 0$ and $C > 0$ so that the following holds: Denote by $H(f_{n,i}) \subset \mathcal{B}_{nor}(U)$ the subset of maps $f \in \mathcal{B}_{nor}(U)$ which are hybrid conjugated with $f_{n,i}$ (defined in proposition 15.1). Then, for $n$ large enough, each connected component of the intersection of $H(f_{n,i})$ with*

$$Q_\epsilon := \{\tilde{f}_i + v_1 + v_2\colon\ v_1 \in E^h_{\tilde{f}_i},\ v_2 \in E^v_{\tilde{f}_i},\ |v_1|_{\mathcal{B}(U)} \leq \epsilon \text{ and } |v_2|_{\mathcal{B}(U)} \leq C\epsilon\}$$

*is a connected complex analytic manifold which is the graph of a complex analytic function*

$$h_{n,i} \colon \{v_1 \colon v_1 \in E^h_{\tilde{f}_i}, \ |v_1| \le \epsilon\} \to \{v_2 \colon v_2 \in E^v_{\tilde{f}_i}, \ |v_2| \le C\epsilon\}.$$

*Note that $\epsilon$ does not depend on $n$.*

The following corollary is a monotonicity result:

**Corollary 16.1.** *$H(f_{n,i}) \cap Q_\epsilon$, for $n$ large enough, has an unique connected component.*

*Proof.* If $c_\infty$ is the Fibonacci parameter, then the family $x^d + c$ is transversal with $E^h_{x^d+c_\infty}$ (because this family is tangent to the vertical direction). By the contraction of the Fibonacci renormalization operator, there exists an induced transformation $\mathcal{I}$ from a neighborhood of $x^d + c_\infty$ in some Banach slice $\mathcal{B}_{nor}(V)$ ($V$ is simply connected) to $\mathcal{B}_{nor}(U)$ so that $\mathcal{I}(x^d + c_\infty) \in W^s(\tilde{f}_1)$ and $D\mathcal{I}(E^h_{x^d+c_\infty}) \subset E^h_{\mathcal{I}(x^d+c_\infty)}$. Furthermore, by the infinitesimal pullback argument and the non existence of invariant line fields on $K(x^d+c_\infty)$, it is easy to see that $D\mathcal{I}^{-1}(E^h_{\mathcal{I}(x^d+c_\infty)}) = E^h_{x^d+c_\infty}$. So the family $\mathcal{I}(x^d + c)$ is transversal to $W^s_{loc}(\tilde{f}_1)$. In particular, if $\epsilon$ is small enough and $n$ is large, the real family $\mathcal{I}(x^d + c)$ intercept each connected component of $H(\tilde{f}_{n,i}) \cap Q_\epsilon$ at least once. But this real family cannot hit $H(f_{n,i})$ in two parameters $c_1$, $c_2$, since $x^d + c_1$ and $x^d + c_2$ will have the same combinatorics and, by a well known result, for each unimodal combinatorics where the critical point is periodic, there exists only one real parameter $c$ so that $x^d + c$ (recall that $d$ is an even number) have this combinatorics. □

*Proof of Theorem 3.* Let $H(f_{n,i})$; $i = 1, 2$; be the codimension one complex manifold of functions which are hybrid conjugated to $f_{n,i}$. Note that

$$\phi(H(f_{n,i}) \cap B_i) = H(f_{n,3-i}) \cap B_{3-i},$$

and

$$\mathcal{R}^{-1}(H(f_{n-1,i}) \cap B_i) \cap B_{3-i} \subset H(f_{n,3-i}),$$

so

$$(\mathcal{R} \circ \phi)^{-1}(H(f_{n-1,i}) \cap B_i) \cap B_i \subset H(f_{n,i}).$$

Let $f_\lambda \in \mathcal{B}_{nor}(U_0 \cup U_1)$, $\lambda \in \mathbb{C}$, be a complex analytic family so that $f_\lambda(z)$ is real if $\lambda$ and $z$ is real. Assume that $f_0$ has the Fibonacci combinatorics. By the topological convergence of the Fibonacci renormalization operator, we can define an induced transformation $\mathcal{I}$ of the neighborhood of $f_0$ in a neighborhood of a point in the local stable manifold of $\tilde{f}_1$. Since the image of the derivative of any induced transformation is dense, we can find a vector $v \in T\mathcal{B}_{nor}(U_0 \cup U_1)$ so that the family $\mathcal{I}(f_\lambda + \lambda \cdot v)$ is transverse to the $W^s_{loc}(\tilde{f}_1)$.

The tangent space of $H(f_{n,i})$ in a point $g$ is exactly $E^h_g$. So $H(f_{n,i})$ is transversal with the unstable direction of $\tilde{f}_i$. Now apply to the operator $\mathcal{R} \circ \phi$ and its fixed points $\tilde{f}_i$ a well known fact about (complex) hyperbolic fixed points with (complex) one dimensional complex unstable manifold, which implies that, if $f_\lambda$ is a curve in $\mathcal{B}_{nor}(U)$ that intersect $W^s_{loc}(\tilde{f}_i)$ in a transversal way at $\lambda = 0$, then there exists, for $n$ large enough, an unique sequence $\lambda_n$ (the unicity follows of corollary 16.1) so that $f_{\lambda_n} \in H(f_{n,i})$ and

$$\frac{\lambda_n - \lambda_{n+1}}{\lambda_{n+1} - \lambda_{n+2}}$$

converges to the unique eigenvalue of $D(\mathcal{R}_{\tilde{f}} \circ \phi)$ larger than one. Because $(\mathcal{R} \circ \phi)^2 = \mathcal{R}^2$, the absolute value of this eigenvalue is exactly the square root of the modulus of the unique eigenvalue larger than one of the operators $D\mathcal{R}^2_{\tilde{f}_1}$ and $D\mathcal{R}^2_{\tilde{f}_2}$ (these two linear operators have the same spectrum, once $\mathcal{R}^2 \circ \phi = \phi \circ \mathcal{R}^2$).

For the family of generalized polynomials induced by the family $x^d + c$, it is not necessary a perturbation: consider the induced transformation defined in the proof of corollary 16.1. We concluded there that the family $\mathcal{I}(x^d + c)$ is transversal to $W^s_{loc}(\tilde{f}_1)$. This finished the proof. $\square$

## 17. Further developments and open questions

The methods explained in this work give some light to

**Conjecture 1.** *The set of parameters $c \in \mathbb{R}$ so that $x^d + c$ has a wild attractor has zero Lebesgue measure.*

H. Bruin[Bru] described a lot of non renormalizable unimodal combinatorics which admits an wild attractor to large criticalities. It seems that in all these examples, the post critical set has bounded geometry, but, as it was noted by Bruin, the inclusion of almost parabolic points does not seems to have any effect in the existence of wild attractors. If a combinatorics has (essentially) bounded geometry (in the W. Shen sense[Sh]) for large criticalities, then it also admits a wild attractor for (perhaps larger) criticalities? (if a non renormalizable map has bounded geometry then it is not necessarily true that this map has a wild attractor: Keller & Nowicki[KN] proved that there exist Fibonacci maps, with criticality $2 + \epsilon$ and bounded geometry, but with an absolutely continuous invariant probability) We expect an affirmative answer. The next step is to see if these examples exhaust the combinatorics which admits wild attractors: Note that wild attractors cannot coexist with linear decay of geometry [Lyu94]. If a map have 'large bounds' for infinitely many levels in the principal nest (exactly the opposite situation of "essentially bounded geometry" in W. Shen work[Sh]), this map can have a wild attractor? A negative answer seems to be so good to be true. Computer experiments and partial results by the author[Sm] suggest that the complex bounds and puzzle geometry control in Theorem 1 holds for all non renormalizable unimodal combinatorics which admits a wild attractor (for large criticalities) described in Bruin[Bru], at least in criticalities where the decay of geometry or bounded geometry holds (certainly for degree two). We expect that all combinatorics which admits essentially bounded geometry for some criticality satisfies the complex bounds and puzzle geometry control in Theorem 1 (possibly not just to the criticalities where the essentially bounded geometry holds, but to all criticalities). So we expect a large intersection (coincidence?) of

- The combinatorial types which satisfies the complex bounds and puzzle geometry control in Theorem 1.
- The combinatorial types which have wild attractors for some criticality.
- The combinatorial types which have (essentially) bounded geometry for some criticality.

Once we have such (partial) identification, the idea is to prove the existence of a hyperbolic set to the generalized renormalization operator associated with these combinatorics (to include combinatorics with almost parabolic fixed points, something similar to the Lyubich "full horseshoe" is necessary), so that this set contains all limit behaviours of the action of the generalized renormalization operator on the maps which are infinitely renormalizable with respect to it. In view of the Lyubich work in the quadratic family [Lyu99], the existence of hyperbolic sets for generalized renormalization operators implies that the parameters in the real family $x^d + c$ whose correspondent map is infinitely renormalizable with respect to these generalized renormalization operators have zero Lebesgue measure in the parameter space (in the families $x^d + c$). The methods applied in this work to the Fibonacci renormalization operator can be generalized to a lot of other combinatorics, including non periodic ones. Note that, to prove the hyperbolicity, it was necessary more than the puzzle geometry control: it is necessary that after normalize the domains of the generalized renormalizations in an appropriated way, we obtain polynomial-like maps which are hyperbolic, in such way we can control the critical puzzle geometry of maps whose renormalization are very close to the original map. We also expect that this kind of argument can be generalized to many combinatorics, incluing non periodic ones.

We hope to give a first approximation to Conjecture 1 in [Sm], using the approach suggested above to prove that sets of parameters corresponding to large sets of combinatorics which admits wild attractors have zero Lebesgue measure in the parameter space.

Another future development involves rigidity results in the real family $x^d + c$, which is related with, for instance, the density of the hyperbolic maps in this real family. For non renormalizable quadratic polynomials, the decay of geometry holds. But for higher criticalities, this is not true: specially bad examples are non renormalizable real maps for which neither decay of geometry nor (essentially) bounded geometry holds: the previous methods to prove rigidity apparently can not be used in this situation. Our Rigidity Theorem will be probably very useful here: for instance, assume that $f$ is a non renormalizable map of even criticality larger than two. Furthermore, assume that $f$ has, in the principal nest,

- A infinity number of finite sequences of generalized renormalizations (cascades) where the decay of geometry does not holds (as in the Fibonacci cascades),
- A infinity number of finite cascades for which decay of geometry holds.

Depending how we combine these two types of cascades, we will get different metric behaviors (decay of geometry, bounded geometry, or any of them). Indeed, we conjecture that we can deal with all cases in the same way: suppose by induction that we have controlled the puzzle geometry of the principal nest at the nth cascade. If in the next cascade the decay of geometry holds, use the Thurston map as in Lyubich work to control the geometry in it. Otherwise, use the Lyubich argument [Lyu93] until the modulus in the principal nest down , and then use arguments as Buff [Bu] (indeed, we can simplify this part). Then use theorem 1 to prove rigidity for these combinatorics. Certainly it is possible to construct by hand a map with neither decay of geometry nor bounded geometry, and so that we can apply the above argument, but we do not know how general the method actually is.

## Acknowledgments

I would like to thank the hospitality of IMS-SUNY-Stony Brook, ICTP-Trieste, Mathematics Institute of Warwick University, and IMPA, where parts of this work were done. I would like to thank specially to M. Lyubich, for the my exciting period at Stony Brook, and S. van Strien, for my wonderful stay at Warwick. Thanks to A. Avila for discussions about Avila, Lyubich and de Melo work and the correction of a mistake in the proof of lemma 12.1. I am also grateful to Juan Rivera-Letelier, whose comments were useful to improve the appresentation of this work, by the fruitful suggestion of an inductive argument with the hyperbolic metric in the proof of Lemma 11.1, and by the great support in the Fall/2001.

Institute for Mathematical Sciences, State University of New York at Stony Brook, Stony Brook-NY, 11794

*E-mail address*: `smania@math.sunysb.edu`